\documentclass{gtart_h}  

\def\ifplaintex{\expandafter\ifx\csname documentclass\endcsname\relax}

\ifplaintex 
\else
\fi

\expandafter\ifx\csname beginpicture\endcsname\relax
\expandafter\ifx\csname documentclass\endcsname\relax
\input pictex \else
\input prepictex \input pictex \input postpictex \fi\fi

\def\gt{{\mathsurround=0pt\it $\cal G\mskip-2mu$eometry \&\ 
$\cal T\!\!$opology}}        

\def\gtp{{\mathsurround=0pt\it $\cal G\mskip-2mu$eometry \&\ 
$\cal T\!\!$opology $\cal P\!$ublications}}  

\def\lognumber#1{\def\thelognumber{#1}}
\def\volumenumber#1{\def\thevolumenumber{#1}}
\def\papernumber#1{\def\thepapernumber{#1}}
\def\volumeyear#1{\def\thevolumeyear{#1}}

\def\pagenumbers#1#2{\def\startpage{#1}\def\finishpage{#2}}
\def\published#1{\def\publishdate{#1}}
\def\proposed#1{\def\theproposer{#1}}
\def\seconded#1{\def\theseconders{#1}}
\def\received#1{\def\receiveddate{#1}}
\def\revised#1{\def\reviseddate{#1}}
\def\accepted#1{\def\accepteddate{#1}}

\long\def\asciiabstract#1{\long\def\theasciiabstract{#1}}
\def\asciikeywords#1{\def\theasciikeywords{#1}}

\def\shorttitle#1{\def\theshorttitle{#1}}

\let\\\par\let\thelognumber\relax
\let\thevolumenumber\relax\let\thepapernumber\relax
\let\thevolumeyear\relax\let\thesamplenumber\relax\let\startpage\relax
\let\finishpage\relax\let\publishdate\relax\let\receiveddate\relax
\let\reviseddate\relax\let\accepteddate\relax\let\theasciititle\relax
\let\theasciiauthors\relax
\let\theasciiabstract\relax\let\theasciikeywords\relax
\let\theasciiemail\relax\let\theshortauthors\relax\let\theshorttitle\relax

\long\def\maketitlep{   

\count0=\startpage

\gt\hfill      
\beginpicture
\setcoordinatesystem units <0.33truein, 0.33truein> point at 2.2 0.9
\setplotsymbol ({$\cal G$})
\plotsymbolspacing=9truept
\circulararc 315 degrees from 0 1 center at 0 0
\setplotsymbol ({$\cal T$})
\circulararc 315 degrees from 1 -1 center at 1 0
\endpicture
\break
{\small\ifx\thesamplenumber\relax 
Volume \else Sample
\fi\thevolumenumber\ (\thevolumeyear)
\startpage--\finishpage\nl
Published: \publishdate}
\vglue 0.5truein plus 0.4fil minus 0.1truein

{\parskip=0pt\leftskip 0pt plus 1fil\def\\{\par\smallskip}{\ifplaintex\large
\else\Large\fi\bf\thetitle}\par\medskip}   

\vglue 0pt plus 0.1fil 

{\parskip=0pt\leftskip 0pt plus 1fil\def\\{\par}{\sc\theauthors}
\par\medskip}

\vglue 0pt plus 0.1fil 

{\small\parskip=0pt\let\newline\\
{\leftskip 0pt plus 1fil\def\\{\par}{\sl\theaddress}\par}
\expandafter\ifx\theemail\relax    
\relax\else\vglue 5pt plus 0.02fil minus 2pt\def\\{\stdspace{\rm 
and}\stdspace} 
\cl{Email:\stdspace\tt\theemail}\fi
\ifx\theurl\relax                  
\relax\else\vglue 5pt plus 0.02fil minus 2pt\def\\{\stdspace{\rm 
and}\stdspace}
\cl{URL:\stdspace\tt\theurl}\fi\par}

\vglue 7pt plus 0.3fil minus 3pt

{\bf Abstract}
\vglue 5pt plus 0.1fil minus 2pt

\theabstract

\vglue 7pt plus 0.3fil minus 3pt

{\bf AMS Classification numbers}\quad Primary:\quad \theprimaryclass

Secondary:\quad \thesecondaryclass

\vglue 5pt plus 0.3fil minus 2pt

{\bf Keywords:}\quad \thekeywords

\vglue 10pt plus 0.5fil minus 5pt

{\small  Proposed: \theproposer\hfill Received: \receiveddate\nl
Seconded: \theseconders\hfill 
\ifx\reviseddate\relax                         
Accepted: \accepteddate                        
\else
Revised: \reviseddate                          
\fi}
\eject
}       

\let\maketitlepage\maketitlep
\let\maketitle\maketitlepage

\font\phead=cmsl9 scaled 950
\font=cmsl9 scaled 1050
\font\pnum=cmbx10 scaled 913
\font\lnum=cmbx10 
\font\pfoot=cmsl9 scaled 950
\font=cmsl9 scaled 1050
\ifplaintex
\headline{\vbox to 0pt{\vskip -4.5mm\line{\small\phead\ifnum
\count0=\startpage ISSN 1364-0380 (on line)
1465-3060 (printed) \hfill {\pnum\folio}\else\ifodd\count0\def\\{ }%
\ifx\theshorttitle\relax\thetitle\else\theshorttitle\fi\hfill{\pnum\folio}
\else\def\\{ and }{\pnum\folio}\hfill\ifx\theshortauthors\relax\theauthors
\else\theshortauthors\fi\fi\fi}\vss}}
\footline{\vbox to 0pt{\vglue 0mm\line{\small\pfoot\ifnum\count0=\startpage
\copyright\ \gtp\hfill\else
\gt, Volume \thevolumenumber\ (\thevolumeyear)\hfill\fi}\vss
}}
\else
\makeatletter
\def\@oddhead{{\small\lhead\ifnum\count0=\startpage ISSN 1364-0380 (on line)
1465-3060 (printed) \hfill {\lnum\number\count0}\else\ifodd\count0
\def\\{ }\ifx\theshorttitle\relax \thetitle \else\theshorttitle\fi\hfill
{\lnum\number\count0}\else\def\\{ and }{\lnum\number\count0}
\hfill\ifx\theshortauthors\relax 
\theauthors\else\theshortauthors\fi\fi\fi}}\def\@evenhead{\@oddhead}
\def\@oddfoot{\small\lfoot\ifnum\count0=\startpage\copyright\ \gtp\hfill\else
\gt, Volume \thevolumenumber\ (\thevolumeyear)\hfill\fi}
\def\@evenfoot{\@oddfoot}
\makeatother
\fi

\newwrite\gtoutfile
\long\gdef\makeheadfile{  
{\def\\{, }\def\s{ }
\immediate\openout\gtoutfile head.xxx
\immediate\write\gtoutfile{Proxy-for: \ifx\theasciiauthors\relax
\theauthors\else\theasciiauthors\fi\s<\ifx\theasciiemail\relax\theemail\else\theasciiemail\fi>}
\immediate\write\gtoutfile{\noexpand\\}
\immediate\write\gtoutfile{Authors: \ifx\theasciiauthors\relax
\theauthors\else\theasciiauthors\fi}
{\def\\{ }\immediate\write\gtoutfile{Title: \ifx\theasciititle\relax
\thetitle\else\theasciititle\fi}}
\immediate\write\gtoutfile{Subj-class: GT or SG or MG etc}
\immediate\write\gtoutfile{MSC-class: \theprimaryclass\ifx\thesecondaryclass\relax\else, \thesecondaryclass\fi}
\immediate\write\gtoutfile{Journal-ref: Geom. Topol. \thevolumenumber
(\thevolumeyear) \startpage-\finishpage}
\immediate\write\gtoutfile{Comments: Published by Geometry and Topology at}
\immediate\write\gtoutfile{\s\s http://www.maths.warwick.ac.uk/gt/GTVol\thevolumenumber/paper\thepapernumber.abs.html}
\immediate\write\gtoutfile{\noexpand\\}
\immediate\write\gtoutfile{}
\ifx\theasciiabstract\relax
\immediate\write\gtoutfile{\theabstract}\else
\immediate\write\gtoutfile{\theasciiabstract}\fi
\immediate\write\gtoutfile{}
\immediate\write\gtoutfile{\noexpand\\}
\immediate\write\gtoutfile{}
\immediate\closeout\gtoutfile}}  

\def\maketitlepage{\maketitlep\makeheadfile}
\let\maketitle\maketitlepage

\lognumber{402}
\volumenumber{8}\papernumber{18}\volumeyear{2004}
\pagenumbers{701}{734}
\received{27 December 2003}
\accepted{3 May 2004}
\revised{20 January 2004}
\published{16 May 2004}
\proposed{Yasha Eliashberg}
\seconded{Robion Kirby, Ronald Fintushel}

\usepackage{amssymb, amsmath}

\newtheorem{thm}{Theorem}[section]    
\newtheorem{lem}[thm]{Lemma}          
        
\newtheorem{coro}[thm]{Corollary}        

\theoremstyle{definition}

\newtheorem{re}[thm]{Remark}

\newcommand{\C}{{\mathbb C}}
\newcommand{\R}{{\mathbb R}}
\newcommand{\Z}{{\mathbb Z}}
\newcommand{\Q}{{\mathbb Q}}
\renewcommand{\P}{{\mathbb P}}
\newcommand{\s}{{\mathbb S}}
\newcommand{\B}{{\mathbb B}}

\newcommand{\M}{{\cal M}}

\newcommand{\E}{{\cal E}}

\renewcommand{\O}{{\cal O}}

\renewcommand{\H}{{\cal H}}

\def\sh#1{\subsection*{#1}\addcontentsline{toc}{subsection}{#1}}

\begin{document}
\title{Orbifold adjunction formula and symplectic\\cobordisms 
between lens spaces}
\shorttitle{Orbifold adjunction formula and symplectic cobordisms}

\author{Weimin Chen}
\address{Mathematics Department, Tulane University\\New Orleans, 
LA 70118, USA}
\email{wchen@math.tulane.edu}

\begin{abstract}
Each lens space has a canonical contact structure which lifts to the
distribution of complex lines on the three-sphere. In this paper, we
show that a symplectic homology cobordism between two lens spaces,
which is given with the canonical contact structure on the boundary,
must be diffeomorphic to the product of a lens space with the unit
interval. As one of the main ingredients in the proof, we also derive
in this paper the adjunction and intersection formulae for
pseudoholomorphic curves in an almost complex 4--orbifold, extending
the relevant work of Gromov and McDuff in the manifold setting.
\end{abstract}

\asciiabstract{%
Each lens space has a canonical contact structure which lifts to the
distribution of complex lines on the three-sphere. In this paper, we
show that a symplectic homology cobordism between two lens spaces,
which is given with the canonical contact structure on the boundary,
must be diffeomorphic to the product of a lens space with the unit
interval. As one of the main ingredients in the proof, we also derive
in this paper the adjunction and intersection formulae for
pseudoholomorphic curves in an almost complex 4-orbifold, extending
the relevant work of Gromov and McDuff in the manifold setting.}

\keywords{Cobordism of lens spaces, orbifold adjunction formula,
symplectic 4--orbifolds, pseudoholomorphic curves}

\asciikeywords{Cobordism of lens spaces, orbifold adjunction formula,
symplectic 4-orbifolds, pseudoholomorphic curves}

\primaryclass{57R17}
\secondaryclass{57R80}

\maketitlepage

\section{Introduction}

In this paper, we prove the following theorem.

\begin{thm}
Let $(W,\omega)$ be a symplectic homology cobordism between two lens spaces
which are equipped with their canonical contact structure. Then $W$ is
diffeomorphic to the product of a lens space with the unit interval.
\end{thm}

Here the canonical contact structure $\xi_0$ on a lens space $L(p,q)$ is the
descendant of the distribution of complex lines on $\s^3=\{(z_1,z_2)\mid
|z_1|^2+|z_2|^2=1\}$ under the quotient map $\s^3\rightarrow L(p,q)$ of the
$\Z_p$--action $(z_1,z_2)\mapsto (\mu_p z_1,\mu_p^q z_2)$. (Here $\mu_p=
\exp(\sqrt{-1}\frac{2\pi}{p})$, and $p,q$ are relatively prime and $0<q<p$.)
The contact structure $\xi_0$ induces a canonical orientation on $L(p,q)$
where a volume form is given by $\alpha\wedge d\alpha$ for some 1--form
$\alpha$ such that $\xi_0=\ker\alpha$. A symplectic cobordism from
$(L(p^\prime,q^\prime),\xi_0^\prime)$ to $(L(p,q),\xi_0)$ is a symplectic
4--manifold $(W,\omega)$ with boundary $\partial W=L(p,q)
-L(p^\prime,q^\prime)$, such that there exists a vectorfield $v$ in
a neighborhood of $L(p,q)\cup L(p^\prime,q^\prime)\subset W$,
which is transverse to $L(p,q)\cup L(p^\prime,q^\prime)$ and for which
$L_v\omega=\omega$, $\xi_0^\prime=\ker\;(i_v\omega|_{L(p^\prime,q^\prime)})$,
$\xi_0=\ker\;(i_v\omega|_{L(p,q)})$, and the canonical orientations on
$L(p,q), L(p^\prime,q^\prime)$ agree with the orientations
defined by the normal vector $v$. (Here $W$ is canonically oriented by
the symplectic form $\omega$, ie, $\omega\wedge\omega$ is a volume form.)
The cobordism $W$ is called a homology cobordism if each $L(p,q)\subset W,
L(p^\prime,q^\prime)\subset W$ induces an isomorphism on homology groups
(with $\Z$ coefficients). In particular, this condition implies $p=p^\prime$.

As a special case, consider the following:

\begin{coro}
Let $\rho$ be a symplectic $\Z_p$--action on $(\R^4,\omega_0)$ where
$\omega_0=dx_1\wedge dy_1+dx_2\wedge dy_2$. Suppose outside of a ball,
$\rho$ is linear and free, and is orthogonal with respect to the Euclidean
metric $g_0=\sum_{i=1}^2 (dx_i^2+dy_i^2)$. Then $\rho$ is conjugate to a 
linear action by a diffeomorphism which is identity outside of a ball. 
\end{coro}

\begin{re}

(1)\qua It is likely that Corollary 1.2 can be strengthened to the assertion
that the action $\rho$ is conjugate to a linear action by a
symplectomorphism of $(\R^4,\omega_0)$. We plan to address this problem
in a separate paper.

(2)\qua Relevant to Theorem 1.1 and Corollary 1.2, we mention two earlier
results. One is due to Eliashberg (cf \cite{E}) which says that a
symplectic $4$--manifold $W$ with contact boundary $\s^3$ (in the weak sense)
is diffeomorphic to a blowup of the $4$--ball $\B^4$. The other is due to
Gromov--McDuff (cf for example Theorem 9.4.2 in \cite{McDS}) which says that if
$(W,\omega)$ is a minimal symplectic $4$--manifold and there are
compact subsets $K\subset W$ and $V\subset \R^4$ with $V$ being star-shaped
such that $(W\setminus K,\omega)$ is symplectormorphic to
$(\R^4\setminus V,\omega_0)$ via a map $\psi$, then there exists a
symplectomorphism $\phi\co (W,\omega)\rightarrow (\R^4,\omega_0)$ which
agrees with $\psi$ on $W\setminus K^\prime$ for some larger compact
subset $K^\prime\supset K$. 

(3)\qua Symplectic fillings (in the weak sense) of lens spaces with
the canonical contact structure are classified up to
orientation-preserving diffeomorphisms in \cite{Lis}, where it is
shown that there are infinitely many lens spaces which have a unique
filling up to blowups. For these lens spaces, it is clear that when
the two ends of $\partial W$ are diffeomorphic, the condition that $W$
is a homology cobordism is equivalent to the condition that
$(W,\omega)$ is minimal.
\end{re}

The proof of Theorem 1.1 is based on studying pseudoholomorphic curves
in a certain symplectic 4--orbifold in the fashion of Gromov--McDuff in
the manifold setting (cf for example \cite{McDS}). There are two main ingredients.
One is the orbifold analog of the adjunction and intersection formulae for
pseudoholomorphic curves, extending the relevant work of Gromov and McDuff
\cite{Gr, McD1, McD2} in the manifold setting. The other is a structural
theorem for the space of a certain notion of maps\footnote{A prototype of
this notion appeared first in \cite{CR} in the disguise 
of ``a good (V--manifold)
map + an isomorphism class of the pull-backs of the tangent bundle''.}
between orbifolds developed in \cite{C1}, which is needed here for the
corresponding Fredholm theory.

The paper is organized as follows. In Section 2 we introduce a notion of
differentiable chains in orbifolds, which serves as a bridge between
the de Rham cohomology of an orbifold and the singular cohomology of
its underlying space via integration. Section 3 is devoted to the proof
of the orbifold analog of the adjunction and intersection formulae.
The main results are proved in Section~4.

\rk{Acknowledgments}

I am indebted to Ron Fintushel for bringing the problem of
$\s^1$--actions on $\s^5$ to my attention, and to S\l awomir Kwasik
for pointing out \cite{Pa} to me, which eventually led to the study in
this paper. I am also very grateful to Dusa McDuff for helpful
communications regarding her relevant work and for kindly letting me
use \cite{McDS} before its publication, and to Morris Kalka and S\l
awomir Kwasik for very valuable conversations. An earlier version of
this work contains a serious mistake, thanks to Slava Matveyev for
pointing it out to me. Finally, I wish to thank an anonymous referee
for pointing out several misspellings in the text and whose queries
helped improve the presentation of this article. This research was
partially supported by NSF Grant DMS-0304956.

\section{Differentiable chains in orbifolds}

We introduce here a notion of differentiable chains in orbifolds. The homology
groups of the corresponding chain complex are naturally isomorphic to the
singular homology groups of the underlying space over $\Q$, so that this
construction yields an explicit pairing between the de Rham cohomology groups
of the orbifold and the singular homology groups of the underlying space via
integration over differentiable chains. In light of the development in
\cite{C1}, the notion introduced here may be regarded as a natural
generalization to the orbifold category of the notion of differentiable 
singular chains in smooth manifolds.

A differentiable $r$--chain in an orbifold $X$ (of class $C^l$ for some
$l\geq 1$) is a finite linear combination of differentiable $r$--simplexes
in $X$, where a differentiable $r$--simplex $\sigma$ in $X$ is a
differentiable map (in the sense of \cite{C1}) from a certain $r$--dimensional
orbihedron into $X$. More precisely, the said $r$--dimensional orbihedron is
an orbispace where the underlying space is the standard $r$--simplex
$\Delta^r$ in $\R^r$, and the orbispace structure is given by a complex
of finite groups over $\Delta^r$ in the sense of Haefliger \cite{Hae}
(see also Part II of \cite{C1}). Recall that a complex of groups consists 
of the
following data: $(K,G_\tau,\psi_a,g_{a,b})$, where $K$ is a simplicial
complex, $G_\tau$ is a group assigned to each cell $\tau\in K$,
$\psi_a\co G_{i(a)}\rightarrow G_{t(a)}$ is an injective homomorphism assigned
to each edge $a$ in the barycentric subdivision of $K$ with $i(a)$, $t(a)$
being the cells of $K$ whose barycenters are the end points of $a$ such that
$t(a)$ is a face of $i(a)$, and $g_{a,b}$ is an element of $G_{t(a)}$ assigned
to each pair of composable edges $a,b$ such that
$$
Ad(g_{a,b})\circ \psi_{ab}=\psi_a\circ\psi_b, \hspace{2mm}
\psi_a(g_{b,c})g_{a,bc}=g_{a,b}g_{ab,c}.
$$
The orbihedron is covered by a set of ``uniformizing systems'' which
are given with compatible equivariant simplicial structures. The $r$--simplex
$\sigma$ being a differentiable map means that the representatives of
$\sigma$ are differentiable when restricted to each simplex in the
corresponding uniformizing system.

Let $\alpha$ be a differential $r$--form on $X$. Then a differentiable
$r$--simplex $\sigma$ in $X$ pulls back $\alpha$ to a differential $r$--form
$\sigma^\ast\alpha$ on $\Delta^r$, the standard $r$--simplex in $\R^r$.
We define the integration of $\alpha$ over $\sigma$ by
$$
\int_\sigma\alpha=\frac{1}{|G|}\int_{\Delta^r}\sigma^\ast\alpha
$$
where $|G|$ is the order of the group $G$ assigned to the top cell
of $\Delta^r$ in the complex of finite groups that defines the orbispace
structure of the orbihedron over which $\sigma$ is defined.
The integration over a differentiable $r$--chain $c=\sum_k a_k\sigma_k$
is defined to be
$$
\int_c\alpha=\sum_k a_k\int_{\sigma_k}\alpha.
$$

Next we introduce a boundary operator $\partial$ on the set of
differentiable chains. To this end, let $\Delta_i^r$, $0\leq i\leq r$,
be the $i$-th face of the standard $r$--simplex $\Delta^r$. The restriction
of a differentiable $r$--simplex $\sigma$ to $\Delta_i^r$ (given the
suborbihedron structure, cf \cite{C1}) is a differentiable $(r-1)$--simplex,
which will be denoted by $\sigma_i$. We define
$$
\partial \sigma=\sum_{i=0}^r (-1)^i\frac{|G_i|}{|G|}\sigma_i
$$
where $G_i,G$ are the groups assigned to the top cell of
$\Delta_i^r,\Delta^r$ respectively. The boundary of a differentiable
$r$--chain $c=\sum_k a_k\sigma_k$ is defined
to be $\partial c=\sum_k a_k\partial\sigma_k$, which clearly satisfies
$$
\partial\circ\partial = 0.
$$
Finally, the Stokes' theorem implies that for any differentiable $r$--chain
$c$ and $(r-1)$--form $\alpha$,
$$
\int_c d\alpha=\int_{\partial c}\alpha.
$$

For any orbifold $X$, let $\H_\ast(X)$, $\H^\ast(X)$ be the homology
and cohomology groups of differentiable chains (with $\Z$ coefficients)
in $X$. There are canonical homomorphisms
$$
H^\ast_{dR}(X)\rightarrow \H^\ast(X)\otimes\R
$$
induced by integration over differentiable chains, and
$$
\H_\ast(X)\rightarrow H_\ast(X;\Q)
$$
which is defined at the chain level by
$$
\sigma\mapsto \frac{1}{|G|}|\sigma|
$$
for each differentiable $r$--simplex $\sigma\co \Delta^r\rightarrow X$, where
$|\sigma|$ is the induced singular $r$--simplex in the underlying space,
and $|G|$ is the order of the group $G$ assigned to the top cell of $\Delta^r$.

\begin{thm}
The canonical homomorphism $H^\ast_{dR}(X)\rightarrow \H^\ast(X)\otimes\R$
is isomorphic, and the canonical homomorphism $\H_\ast(X)\rightarrow
H_\ast(X;\Q)$ is isomorphic over $\Q$.
\end{thm}

Theorem 2.1 will not be used in this paper, and its proof will be given 
elsewhere. But we remark that the key point in the proof is to show that 
$\H^\ast(X)\otimes\Q$ are the cohomology groups associated to a fine 
torsionless resolution of the constant sheaf $\Q\times X$, with which the 
proof follows by the usual sheaf theoretical argument, for instance, as in 
\cite{W}.

In light of Theorem 2.1, we will say that a differentiable cycle $c$ in
$X$ (ie, a differentiable chain $c$ such that $\partial c=0$) is
Poincar\'{e} dual to a de Rham cohomology class $\gamma\in H_{dR}^\ast(X)$
if there is a closed form $\alpha\in\gamma$ such that for any closed form
$\beta$ on $X$,
$$
\int_c \beta=\int_X\alpha\wedge\beta.
$$
Here is a typical situation: Let $Y$ be a compact, closed, and oriented
$r$--dimensional orbifold and $f\co Y\rightarrow X$ be a differentiable map
in the sense of \cite{C1}. Note that $Y$ can be triangulated such that
with respect to the triangulation, $Y$ is natually an orbihedron
(cf Part II of \cite{C1}). Thus the restriction of $f$ to each top 
simplex in the triangulation of $Y$ defines a differentiable $r$--simplex 
in $X$, and in this way $f(Y)$ naturally becomes a differentiable $r$--chain 
in $X$ which is a cycle because $Y$ is compact, closed, and oriented. 
Clearly, in this case we have
$$
\int_{f(Y)}\beta=\int_Y f^\ast\beta
$$
for any differential form $\beta$ on $X$.

\section{Adjunction and intersection formulae}

In this section, we derive the adjunction formula for pseudoholomorphic
curves in an almost complex 4--orbifold and a corresponding formula which
expresses the algebraic intersection number of two distinct pseudoholomorphic
curves in terms of local contributions from their geometric intersection,
extending relevant work of Gromov \cite{Gr} and McDuff \cite{McD1,McD2} in
the manifold setting.

First of all, some convention and terminology. In this section (and the
previous one as well), the notion of orbifolds is more general in the
sense that the group action on each uniformizing system needs not to be
effective. The orbifolds in the classical sense where the group actions
are effective are called reduced. The points which are the principal orbits
in each uniformizing system are called regular points. They have the smallest
isotropy groups in each connected component of the orbifold, which are all
isomorphic, and they form an open, dense submanifold of the orbifold. The
points in the complement of regular points are called orbifold points.
When the orbifold is reduced and has no codimension $2$ subsets of orbifold
points, we also allow ourselves to use the usual terminologies, ie,
``orbifold point'' = ``singular point''
and ``regular point'' = ``smooth point''.

We now begin by setting the stage. Let $X$ be a compact, closed, and almost
complex 4--dimensional orbifold which is canonically oriented by the almost
complex structure $J$. We shall assume that the 4--orbifold $X$ is reduced
throughout. We shall also consider connected, compact, and closed complex
orbifolds $\Sigma$ with $\dim_\C\Sigma=1$, namely the orbifold Riemann
surfaces, which are not assumed to be reduced in general.

\sh{Definition A}

A $J$--holomorphic curve in $X$ is a closed subset $C\subset X$ such that
there is a nonconstant map $f\co \Sigma\rightarrow X$ in the sense of \cite{C1}
with $C=\mbox{Im }f$,\footnote{Each map $f$ in the sense of \cite{C1}
induces a continuous map between the underlying spaces; by the image under
such an $f$, we always mean the image under the map induced by $f$.}
which obeys
\begin{itemize}
\item [{(a)}] The representatives of $f$ are $J$--holomorphic.
\item [{(b)}] The homomorphisms between isotropy groups in each
representative of $f$ are injective, and are isomorphic at all
but at most finitely many regular points of $\Sigma$.
\item [{(c)}] The map $f$ is not multiply covered in the following
sense: $f$ does not factor through any holomorphic map $\phi\co \Sigma
\rightarrow\Sigma^\prime$ to a map $f^\prime\co \Sigma^\prime\rightarrow X$
such that the degree of the map induced by $\phi$ between the underlying
Riemann surfaces is greater than one.
\end{itemize}

A $J$--holomorphic curve $C$ is called of type I if $\Sigma$ is reduced, and
is called of type II otherwise. Clearly this definition is independent of the
parametrization $f\co \Sigma\rightarrow X$. Likewise, the order of the isotropy
groups of the ``regular'' points in $C$, ie, the images of all but
at most finitely many regular points in $\Sigma$ under $f$, depends only
on $C$, and is called the multiplicity of $C$ and is denoted by $m_C$
throughout. A $J$--holomorphic curve $C$ is of type I if and only if $m_C=1$.
A type I $J$--holomorphic curve is contained in the set of regular points
of $X$ except for possibly finitely many points, and a type II $J$--holomorphic
curve is contained entirely in the set of orbifold points of $X$. Finally,
we remark that for a type I $J$--holomorphic curve $C$, any parametrization
$f\co \Sigma\rightarrow X$ of $C$ is uniquely determined by the induced map
between the underlying spaces.

\sh{Definition B}

{(1)}\qua For any $J$--holomorphic curve $C$ in $X$, the Poincar\'{e}
dual of $C$ is defined to be the class $PD(C)\in H^2(X;\Q)$ which is uniquely
determined by
$$
m_C^{-1}\alpha[C]=PD(C)\cup\alpha[X], \forall \alpha\in H^2(X;\Q),
$$
where $[C]$ is the class of $C$ in $H_2(X;\Z)$.

{(2)}\qua The algebraic intersection number of two $J$--holomorphic curves
$C,C^\prime$ (not necessarily distinct) is defined to be
$$
C\cdot C^\prime=PD(C)\cup PD(C^\prime)[X].
$$

We remark that the Poincar\'{e} dual $PD(C)$ differs from the usual one
by a factor $m_C^{-1}$, thus is different for a type II $J$--holomorphic curve.
On the other hand, if $C$ is parametrized by $f\co \Sigma\rightarrow X$,
the class of the differentiable cycle $f(\Sigma)$ in $\H_2(X)$ is sent
to $m_C^{-1}[C]$ under the canonical homomorphism $\H_2(X)\rightarrow
H_2(X;\Q)$. In light of Theorem 2.1, $PD(C)$ is Poincar\'{e} dual to
$f(\Sigma)$ under the canonical isomorphisms $H^2_{dR}(X)\cong \H^2(X)
\otimes\R\cong H^2(X;\R)$.

We proceed further with a digression on some crucial local properties of
$J$--holomorphic curves in $\C^2$ due to McDuff, cf \cite{McD1,McD2},
where we assume that $\C^2$ is given with an almost complex structure
$J$ which equals the standard structure at the origin. To fix the notation,
the disc of radius $R$ in $\C$ centered at $0$ is denoted by $D(R)$.

First, some local analytic properties of $J$--holomorphic curves:
\begin{itemize}
\item For any $J$--holomorphic curve $f\co (D(R),0)\rightarrow
(\C^2,0)$ where $f$ is not multiply covered, there exists an
$0<R^\prime\leq R$ such that $f|_{D(R^\prime)\setminus\{0\}}$ is embedded.
\item Let $f\colon (D(R),0)\!\rightarrow\! (\C^2,0)$ be a $J$--holomorphic
curve such that $f|_{D(R)\setminus\{0\}}$ is embedded. Then for any
sufficiently small $\epsilon>0$, there is an almost complex
structure $J_\epsilon$ and a $J_\epsilon$--holomorphic immersion
$f_\epsilon$ (not multiply covered) such that as $\epsilon\rightarrow 0$,
$J_\epsilon\rightarrow J$ in $C^1$ topology and $f_\epsilon\rightarrow f$
in $C^2$ topology. Moreover, given any annuli $\{\lambda\leq |z|\leq R\}$
and $\{\lambda^\prime\leq |z|\leq\lambda\}$ in $D(R)$, one can arrange
to have $f=f_\epsilon$ in $\{\lambda\leq |z|\leq R\}$ and to
have $J_\epsilon=J$ except in a chosen neighborhood of the image
of $\{\lambda^\prime\leq |z|\leq\lambda\}$ under $f$ by letting
$\epsilon>0$ sufficiently small.
\item Any two distinct $J$--holomorphic curves $f\co D(R)\rightarrow
\C^2$, $f^\prime\co D(R^\prime)\rightarrow \C^2$ intersect at only finitely
many points, ie, the set $\{(z,z^\prime)\in D(R)\times D(R^\prime)\mid
f(z)=f^\prime(z^\prime)\}$ is finite.
\end{itemize}

Second, the local intersection and self-intersection number
of $J$--holomorphic curves:
\begin{itemize}
\item Let $C$, $C^\prime$ be distinct $J$--holomorphic curves which are
parametrized by $f\co (D(R),0)\rightarrow (\C^2,0)$ and $f^\prime\co (D(R^\prime),0)
\rightarrow (\C^2,0)$, such that $f|_{D(R)\setminus\{0\}}$ and
$f^\prime|_{D(R^\prime)\setminus\{0\}}$ are embedded and $0\in \C^2$
is the only intersection of $C$ and $C^\prime$. Perturb $C$ into
$\overline{C}$ (which may not be pseudoholomorphic), keeping 
$\partial C$ and $C$ disjoint from $C^\prime$ and $\partial C^\prime$ 
respectively, such that $\overline{C}$ intersects with $C^\prime$ 
transversely. Then the intersection number $C\cdot C^\prime$ is defined 
by counting the intersection of $\overline{C}$ and $C^\prime$ with signs.
$C\cdot C^\prime$ may be determined using the following recipe: perturb 
$f,f^\prime$ into $J_\epsilon$--holomorphic immersions 
$f_\epsilon,f_\epsilon^\prime$, then
$$
C\cdot C^\prime=\sum_{\{(z,z^\prime)|f_\epsilon(z)
=f_\epsilon^\prime(z^\prime)\}}t_{(z,z^\prime)}
$$
where $t_{(z,z^\prime)}=1$ when $f_\epsilon(z)=f_\epsilon^\prime(z^\prime)$
is a transverse intersection, and $t_{(z,z^\prime)}=n\geq 2$ when
$f_\epsilon(z)=f_\epsilon^\prime(z^\prime)$ has tangency of order
$n$. The intersection number $C\cdot C^\prime$ has the following properties:
it depends only on the germs of $C,C^\prime$ at $0\in\C^2$, it is always
positive, and it equals one if and only if $C,C^\prime$ are both embedded
and intersect at $0\in\C^2$ transversely.
\item Let $C$ be a $J$--holomorphic curve which is parametrized by
$f\co (D(R),0)$ $\rightarrow (\C^2,0)$ such that $f|_{D(R)\setminus\{0\}}$
is embedded. Then the local self-inter\-sec\-tion number $C\cdot C$ is
well-defined, which can be determined using the following recipe:
perturb $f$ into a $J_\epsilon$--holomorphic immersion $f_\epsilon$,
then
$$
C\cdot C=\sum_{\{[z,z^\prime]|z\neq z^\prime,
f_\epsilon(z)=f_\epsilon(z^\prime)\}}t_{[z,z^\prime]},
$$
where $[z,z^\prime]$ denotes the unordered pair of $z,z^\prime$, and
where $t_{[z,z^\prime]}=1$ when $f_\epsilon(z)=f_\epsilon(z^\prime)$
is a transverse intersection, and $t_{[z,z^\prime]}=n\geq 2$ when
$f_\epsilon(z)=f_\epsilon(z^\prime)$ has tangency of order $n$. The local
self-intersection number $C\cdot C$ has the following properties:
it depends only on the germ of $C$ at $0\in\C^2$, and it is non-negative
which equals zero if and only if $C$ is embedded.
\end{itemize}

End of digression.

In order to state the adjunction and intersection formulae, we need to
further introduce some definitions.

(1)\qua Recall from \cite{C1} that a representative
of a map $f\co \Sigma\rightarrow X$ parametrizing a $J$--holomorphic curve
$C$ gives rise to a collection of pairs
$(f_i,\rho_i)\co (\widehat{D_i},G_{D_i})\rightarrow (\widehat{U_i},G_{U_i})$
satisfying certain compatibility conditions, where
$\{(\widehat{D_i},G_{D_i})\}$,\break $\{(\widehat{U_i},G_{U_i})\}$ are a
collection of uniformizing systems of $\Sigma$ and $X$ respectively, and
each $\rho_i$ is a homomorphism, which is injective by (b) of Definition A,
and each $f_i$ is a $\rho_i$--equivariant $J$--holomorphic map. We may assume
without loss of generality that each $\widehat{D_i}$ is a disc centered
at $0\in\C$ and each $\widehat{U_i}$ is a ball centered at $0\in\C^2$,
and $G_{D_i},G_{U_i}$ act linearly. Moreover, because of (b) and (c) in
Definition A, we may assume that each $f_i$ is embedded when restricted
to $\widehat{D_i}\setminus\{0\}$ and $\rho_i(G_{D_i})$ is the subgroup of
$G_{U_i}$ which leaves $f_i(\widehat{D_i})\subset \widehat{U_i}$ invariant.
(The case of type II is explained in the proof of Lemma 3.4 below.)
Let $z$ be the orbit of $0\in\widehat{D_i}$ in $\Sigma$. We shall call the
germ of $\mbox{Im} f_i$ at $0\in\widehat{D_i}$ a local representative of the
$J$--holomorphic curve $C=\mbox{Im} f$ at $z\in\Sigma$. The set $\Lambda(C)_z$
of all local representatives of $C$ at $z$ is clearly the set of germs of the
elements in
$$
\{\mbox{Im}(g\circ f_i)\mid g\in G_{U_i}\},
$$
which is naturally parametrized by the coset $G_{U_i}/\rho_i(G_{D_i})$.
Note that for all but at most finitely many points $z\in\Sigma$, the set
$\Lambda(C)_z$ of local representatives of $C$ at $z$ contains only one
element.

(2)\qua For any $J$--holomorphic curve $C$ in $X$, its virtual genus is
defined to be
$$
g(C)=\frac{1}{2}(C\cdot C+c(C))+\frac{1}{m_C}
$$
where $c=-c_1(TX)$. Note that $g(C)$ is a rational number in
general.

(3)\qua Let $\Sigma$ be an (connected) orbifold Riemann surface, and
let $m_\Sigma$ be the order (of isotropy groups) of its regular points
and $m_1,m_2,\cdots,m_k$ be the orders (of isotropy groups) of its orbifold
points. We define the orbifold genus of $\Sigma$ by
$$
{g}_\Sigma=\frac{g_{|\Sigma|}}{m_\Sigma}+\sum_{i=1}^k(\frac{1}{2m_\Sigma}-
\frac{1}{2m_i}),
$$
where $g_{|\Sigma|}$ is the genus of the underlying Riemann surface of
$\Sigma$. Note that with the above definition, $c_1(T\Sigma)(\Sigma)
=2m_{\Sigma}^{-1}-2g_\Sigma$ where $T\Sigma$ is the orbifold tangent bundle.

\vspace{3mm}

With the preceding understood, consider the following:

\begin{thm}[Adjunction Formula]
Let $C$ be a $J$--holomorphic curve which is parametrized by
$f\co \Sigma\rightarrow X$. Then
$$
g(C)=g_\Sigma + \sum_{\{[z,z^\prime]|z\neq z^\prime,
f(z)=f(z^\prime)\}}k_{[z,z^\prime]} +
\sum_{z\in\Sigma}k_z,
$$
where $[z,z^\prime]$ denotes the unordered pair of $z,z^\prime$, and
where the numbers $k_{[z,z^\prime]},k_z$ are defined as follows.
\begin{itemize}
\item Let $G_{[z,z^\prime]}$ be the isotropy group at $f(z)=f(z^\prime)$
and $\Lambda(C)_z=\{C_{z,\alpha}\}$, $\Lambda(C)_{z^\prime}=
\{C_{z^\prime,\alpha^\prime}\}$, then
$$
k_{[z,z^\prime]}=\frac{1}{|G_{[z,z^\prime]}|}
\sum_{\alpha,\alpha^\prime}C_{z,\alpha}\cdot C_{z^\prime,\alpha^\prime}.
$$
\item Let $G_z$ be the isotropy group at $f(z)$ and $\Lambda(C)_z
=\{C_{z,\alpha}\}$, then
$$
k_z=\frac{1}{2|G_z|}(\sum_\alpha
C_{z,\alpha}\cdot C_{z,\alpha}+\sum_{\alpha,\beta}C_{z,\alpha}\cdot
C_{z,\beta}).
$$
{\em(}Note: the second sum is over all $\alpha,\beta$ which are not
necessarily distinct.{\em)}
\end{itemize}
\end{thm}

\begin{thm}[Intersection Formula]
Let $C,C^\prime$ be distinct $J$--holomorphic curves parametrized by
$f\co \Sigma\rightarrow X$, $f^\prime\co \Sigma^\prime\rightarrow X$ respectively.
Then the algebraic intersection number
$$
C\cdot C^\prime=\sum_{\{(z,z^\prime)|f(z)=f^\prime(z^\prime)\}}
k_{(z,z^\prime)}
$$
where $k_{(z,z^\prime)}$ is defined as follows. Let $G_{(z,z^\prime)}$ be
the isotropy group at $f(z)=f^\prime(z^\prime)$ and $\Lambda(C)_z=
\{C_{z,\alpha}\}$, $\Lambda(C^\prime)_{z^\prime}=
\{C^\prime_{z^\prime,\alpha^\prime}\}$, then
$$
k_{(z,z^\prime)}=\frac{1}{|G_{(z,z^\prime)}|}
\sum_{\alpha,\alpha^\prime}C_{z,\alpha}\cdot
C^\prime_{z^\prime,\alpha^\prime}.
$$
\end{thm}

The adjunction formula implies the following:

\begin{coro}
Let $C$ be a $J$--holomorphic curve parametrized by $f\co \Sigma\rightarrow X$.
Then the virtual genus of $C$ is greater than or equal to the orbifold
genus of $\Sigma$, ie, $g(C)\geq g_\Sigma$, with $g(C)=g_\Sigma$ iff
$C$ is a suborbifold of $X$ and $f$ is an orbifold embedding.
\end{coro}

The rest of this section is occupied by the proof of Theorem 3.1 and
Theorem 3.2. We begin with some preliminary lemmas.

\begin{lem}
Let $C$ be a type II $J$--holomorphic curve parametrized by $f\co \Sigma
\rightarrow X$. Then $f$ is represented by a collection of pairs
$\{(f_i,\rho_i)\}$ where each $f_i$ is an embedding.
\end{lem}

\proof
Let $(\widehat{U},G_{U})$ be a uniformizing system of $X$, where $\widehat{U}$
is a ball in $\C^2$ and $G_{U}$ is nontrivial and acts linearly. We say that
$G_U$ is of type A if the fixed-point set of $G_U$ is a complex line in
$\C^2$, and that $G_U$ is of type B if $0\in\C^2$ is the only fixed point.

Let $\{(f_i,\rho_i)\}$ be a representative of $f$ (cf \cite{C1}), where
each $(f_i,\rho_i)\co (\widehat{D_i},G_{D_i})\rightarrow (\widehat{U_i},
G_{U_i})$. Since $C$ is of type II, each $G_{U_i}$ is nontrivial. Consider
the case where $G_{U_i}$ is of type A first. In this case, $\mbox{Im }f_i$
lies in the complex line which is fixed by $G_{U_i}$, therefore $f_i$
is a holomorphic map between two discs in $\C$. It follows that $f_i$ is
either an embedding or a branched covering. Suppose $f_i$ is a branched
covering, and without loss of generality assume that $0\in\widehat{D_i}$
is the only branching point. Then there are $z,z^\prime\neq 0$ in
$\widehat{D_i}$ with $z\neq z^\prime$, such that $f_i(z)=f_i(z^\prime)\in
\widehat{U_i}$. Since $f$ is not multiply covered, there must be a $g\in
G_{D_i}$ such that $g\cdot z=z^\prime$. On the other hand, by (b) of
Definition A, $\rho_i$ is an isomorphism onto $G_{U_i}$ when restricted to
the isotropy subgroup of $z$, so that there is an $h\in G_{D_i}$ fixing $z$
such that $\rho_i(h)=\rho_i(g)$. It is easily seen that $\rho_i(gh^{-1})=1
\in G_{U_i}$ but $gh^{-1}\neq 1\in G_{D_i}$, a contradiction to the
assumption in (b) of Definition A that $\rho_i$ is injective. Hence $f_i$
is an embedding. When $G_{U_i}$ is of type B, $\mbox{Im} f_i$ lies in a
complex line in $\C^2$ whose isotropy is a proper subgroup $H$ of $G_{U_i}$.
Again $f_i$ is either an embedding or a branched covering. If $f_i$ is
a branched covering, then there are $z,z^\prime\neq 0$ in $\widehat{D_i}$
with $z\neq z^\prime$, such that $f_i(z)=f_i(z^\prime)\in\widehat{U_i}$.
Moreover, since $f$ is not multiply covered, there is a $g\in G_{D_i}$
such that $g\cdot z=z^\prime$, and in this case, note that $\rho_i(g)\in H$.
On the other hand, there is an $h$ in the isotropy subgroup of $z$ such that
$\rho_i(h)=\rho_i(g)\in H$, which gives a contradiction as in the type A
case. Hence the lemma.
\endproof

\begin{lem}
Let $C$ be a $J$--holomorphic curve parametrized by $f\co \Sigma\rightarrow
X$. Then there is a closed $2$--form $\eta_C$ on $X$ which represents the 
Poincar\'{e} dual of the differentiable cycle $f(\Sigma)$ in $X$, ie,
for any $2$--form $\alpha$ on $X$,
$$
\int_\Sigma f^\ast\alpha=\int_X\eta_C\wedge\alpha.
$$
Moreover, $\eta_C$ may be chosen such that it is supported in any 
given neighborhood of $C$ in $X$. 
\end{lem}

\proof
We consider the case where $C$ is of type I first.

To fix the notation, let $z_1,z_2,\cdots,z_k$ be the set of points
in $\Sigma$ whose image under $f$ is an orbifold point in $X$.
For each $i=1,2,\cdots,k$, we set $p_i=f(z_i)$ and let $m_i\geq 1$ be the
order of the isotropy group at $z_i$. Furthermore, we denote by
$(\widehat{D_i},\Z_{m_i})$, $(\widehat{V_i},G_i)$ some local uniformizing
systems at $z_i$, $p_i$ respectively, and denote by $(f_i,\rho_i)\co 
(\widehat{D_i},\Z_{m_i})
\rightarrow (\widehat{V_i},G_i)$ a local representative of $f$ at $z_i$
such that $f_i$ is embedded when restricted to $\widehat{D_i}\setminus\{0\}$. 
Set $D_i=\widehat{D_i}/\Z_{m_i}$ and $V_i=\widehat{V_i}/G_i$ for the 
corresponding neighborhood of $z_i$ and $p_i$ in $\Sigma$ and $X$ 
respectively. Without loss of generality, we may assume that $D_i$ is 
the connected component of $f^{-1}(V_i)$ that contains $z_i$.

For each critical point $z$ of $f$ (ie $df(z)=0$) where $f(z)$
is a regular point in $X$, we perturb $f$ locally in a small neighborhood
of $z$ into a $J_\epsilon$--holomorphic immersion, which is supported in
the complement of $\bigcup_{i=1}^k D_i$, and for each $i=1,2,\cdots,k$,
we perturb $f_i$ into a $J_\epsilon$--holomorphic immersion $f_{i,\epsilon}$
(if $f_i$ is already embedded, we simply let $f_{i,\epsilon}=f_i$).
Let $\widehat{D_i^\prime}\subset\widehat{D_i}$ be a closed disc of
a smaller radius such that $f_{i,\epsilon}=f_i$ over $\widehat{D_i}
\setminus \widehat{D_i^\prime}$. We set $\Sigma_0=\Sigma\setminus
\bigcup_{i=1}^k D_i$ and $\Sigma_0^\prime=\Sigma\setminus\bigcup_{i=1}^k
D_i^\prime$ where $D_i^\prime=\widehat{D_i^\prime}/\Z_{m_i}$, and we denote
the perturbation of $f$ over $\Sigma_0^\prime$ by $f_\epsilon$, which is
a $J_\epsilon$--holomorphic immersion into $X^0$, the complement of orbifold
points in $X$. Note that $f_{i,\epsilon}$ may not be $\rho_i$--equivariant,
and $J_\epsilon$ may not be $G_i$--equivariant over $\widehat{V_i}$. Hence
$f_\epsilon$, $f_{i,\epsilon}$, $i=1,2,\cdots,k$, may not define a
pseudoholomorphic curve in $X$. Nevertheless, for any closed 2--form $\alpha$
on $X$, it is easily seen that
$$
\int_{\Sigma} f^\ast\alpha=\int_{\Sigma_0} f^\ast_\epsilon\alpha+
\sum_{i=1}^k\frac{1}{m_i}\int_{\widehat{D_i}}
f^\ast_{i,\epsilon}\alpha.
$$

Let $\nu_\epsilon=f_\epsilon^\ast TX^0/T\Sigma_0^\prime$ be the normal bundle
of the immersion $f_\epsilon$ in $X^0$, and let $\nu_{i,\epsilon}
=f^\ast_{i,\epsilon}T\widehat{V_i}/T\widehat{D_i}$ be the normal bundle
of the immersion $f_{i,\epsilon}$ in $\widehat{V_i}$, $i=1,2,\cdots,k$.
We fix an immersion $\bar{f}_\epsilon$ of a tubular neighborhood of the
zero section of $\nu_\epsilon$ into $X^0$, and fix an immersion
$\bar{f}_{i,\epsilon}$ of a tubular neighborhood of
the zero section of $\nu_{i,\epsilon}$ into $\widehat{V_i}$ for each $i$,
which are assumed to be compatible on the overlaps. We denote by
$\overline{\Theta_\epsilon}$, $\overline{\Theta_{i,\epsilon}}$ the
push-forward of some Thom forms $\Theta_\epsilon$, $\Theta_{i,\epsilon}$
of $\nu_\epsilon$, $\nu_{i,\epsilon}$ by $\bar{f}_\epsilon$,
$\bar{f}_{i,\epsilon}$ respectively, where $\Theta_\epsilon$,
$\Theta_{i,\epsilon}$ are compatible on the overlaps. Finally, let
$x_1,x_2,\cdots,x_l$ be the set $\{p_i\mid i=1,2,\cdots,k\}$. For each
$x_j$, $j=1,2,\cdots,l$, let $(\widehat{V_{x_j}},G_{x_j})$ be a local
uniformizing system at $x_j$. Without loss of generality, we assume
$V_i=V_{x_j}=\widehat{V_{x_j}}/G_{x_j}$ whenever $p_i=x_j$.

With the preceding understood, the 2--form $\eta_C$ is defined as follows.
On $X\setminus\bigcup_{j=1}^l V_{x_j}$, $\eta_C=\overline{\Theta_\epsilon}$,
and on each $\widehat{V_{x_j}}$, $j=1,2,\cdots,l$,
$$
\eta_C=\sum_{\{i|p_i=x_j\}}\frac{1}{m_i}\sum_{g\in G_{x_j}}
g^\ast\overline{\Theta_{i,\epsilon}}\,.
$$
Now for any 2--form $\alpha$ on $X$, we have
\begin{align*}
\int_X \eta_C\wedge\alpha
& =  \int_{X\setminus\bigcup_{j=1}^l V_{x_j}}\eta_C\wedge\alpha +
\sum_{j=1}^l\frac{1}{|G_{x_j}|}\int_{\widehat{V_{x_j}}}\eta_C\wedge\alpha\\
& =  \int_{X\setminus\bigcup_{j=1}^l V_{x_j}}\overline{\Theta_\epsilon}
\wedge\alpha + \sum_{j=1}^l\frac{1}{|G_{x_j}|}\int_{\widehat{V_{x_j}}}
(\sum_{\{i|p_i=x_j\}}\frac{1}{m_i}\sum_{g\in G_{x_j}}g^\ast
\overline{\Theta_{i,\epsilon}})\wedge\alpha\\
& =  \int_{\Sigma_0}f_\epsilon^\ast\alpha + \sum_{j=1}^l
\sum_{\{i|p_i=x_j\}}\frac{1}{m_i}(\frac{1}{|G_{x_j}|}
\int_{\widehat{V_{x_j}}}\sum_{g\in G_{x_j}}g^\ast
(\overline{\Theta_{i,\epsilon}}\wedge\alpha))\\
& =  \int_{\Sigma_0}f_\epsilon^\ast\alpha + \sum_{j=1}^l
\sum_{\{i|p_i=x_j\}}\frac{1}{m_i}
\int_{\widehat{V_{x_j}}}\overline{\Theta_{i,\epsilon}}\wedge\alpha\\
& =  \int_{\Sigma_0}f^\ast_\epsilon\alpha+ \sum_{i=1}^k\frac{1}{m_i}
\int_{\widehat{D_i}}f^\ast_{i,\epsilon}\alpha = \int_{\Sigma} f^\ast\alpha.
\end{align*}
Hence $\eta_C$ represents the Poincar\'{e} dual of the differentiable
cycle $f(\Sigma)$. By way of construction, $\eta_C$ may be chosen 
to be supported in any given neighborhood of $C$ in $X$.

\vspace{2mm}

Next we consider the case where $C$ is of type II.

By Lemma 3.4, $\nu=f^\ast TX/T\Sigma$ is an orbifold complex
line bundle over $\Sigma$. Let $\Theta$ be a Thom form of $\nu$. Then
notice that $\nu$ is sort of a quasi-normal bundle of $C$ in $X$ in the
sense that one can push-forward $\Theta$ to $X$. The resulting form, which
is defined to be $\eta_C$, is a closed $2$--form on $X$, supported
in any given neighborhood of $C$, and for any $x\in C$, there exists 
a local uniformizing system $(\widehat{V},G)$ at $x$ such that on 
$\widehat{V}$,
$$
\eta_C=\sum_{i=1}^l\frac{1}{m_i}\sum_{g\in G}g^\ast\overline{\Theta_i},
$$
where $f^{-1}(x)=\{z_1,z_2,\cdots,z_l\}$, $m_i$ is the order of $z_i$
in $\Sigma$, and $\overline{\Theta_i}$ is the push-forward of $\Theta$ 
to $\widehat{V}$
associated to some arbitrarily fixed choice of representatives of the
parametrization $f\co \Sigma\rightarrow X$ of $C$. As in the case where $C$
is of type I, we have for any $2$--form $\alpha$ on $X$
$$
\int_X\eta_C\wedge\alpha=\int_\Sigma f^\ast\alpha,
$$
so that $\eta_C$ represents the Poincar\'{e} dual of the differentiable
cycle $f(\Sigma)$.
\endproof

Note that by the above lemma, we have
$$
C\cdot C^\prime=\int_X\eta_C\wedge\eta_{C^\prime}
$$
for the algebraic intersection number of two $J$--holomorphic curves
$C,C^\prime$.

The next lemma is concerned with a formula which expresses
the first Chern class of an orbifold complex vector bundle over a reduced
orbifold Riemann surface in terms of the first Chern class over the
complement of the orbifold points with respect to a certain canonical
trivialization and the ``first Chern class'' at each orbifold point.
To be more precise, let $E\rightarrow\Sigma$ be a rank $n$ orbifold complex
vector bundle over a reduced orbifold Riemann surface. Let $z_1,z_2,\cdots,
z_k\in\Sigma$ be any given set of points which contains the set of orbifold
points, and let $m_1,m_2,\cdots,m_k$ be the orders of the corresponding
isotropy groups. Suppose over a local uniformizing system
$(\widehat{D_i},\Z_{m_i})$ at each $z_i$, the orbifold bundle $E$ has a
trivialization $(\widehat{D_i}\times\C^n,\Z_{m_i})$, such that $\Z_{m_i}$
acts on $\widehat{D_i}\times\C^n$ by
$$
\mu_{m_i}\cdot (z,v_1,v_2,\cdots,v_n)=(\mu_{m_i}z,
\mu_{m_i}^{m_{i,1}}v_1,\mu_{m_i}^{m_{i,2}}v_2,\cdots,
\mu_{m_i}^{m_{i,n}}v_n),
$$
where $\mu_{m_i}=\exp(\sqrt{-1}\frac{2\pi}{m_i})$ is the generator of
$\Z_{m_i}$, and $0\leq m_{i,j}<m_i$, $j=1,2,\cdots,n$. Set
$D_i=\widehat{D_i}/\Z_{m_i}$, $\Sigma_0=\Sigma\setminus\bigcup_{i=1}^k D_i$,
and $E_0=E|_{\Sigma_0}$. We consider the trivialization $\tau$
of $E_0$ over $\partial\Sigma_0=\bigcup_{i=1}^k\partial D_i$ where
along each $\partial D_i$, $\tau$ is given by pushing down a set of
equivariant sections $\{s_j(z)\mid j=1,2,\cdots,n\}$ of
$\partial\widehat{D_i}\times\C^n$ over $\partial\widehat{D_i}$,
where $s_j(z)=(0,\cdots,z^{m_{i,j}},\cdots,0)$, $j=1,2,\cdots,n$.
Let $\partial D_i\times\C^n$ be the trivialization $\tau$ of $E_0$
over $\partial D_i$. Then the canonical map $\psi_i\co \partial\widehat{D_i}
\times\C^n\rightarrow \partial D_i\times\C^n$ is given by
$$
\psi_i(z,v_1,v_2,\cdots,v_n)=(z^{m_i},z^{-m_{i,1}}v_1,z^{-m_{i,2}}v_2,
\cdots,z^{-m_{i,n}}v_n).
$$
With the preceding understood, the said formula is the following:

\begin{lem}
$c_1(E)(\Sigma)=c_1(E_0,\tau)(\Sigma_0,\partial\Sigma_0)+\sum_{i=1}^k
(\sum_{j=1}^n\frac{m_{i,j}}{m_i})$.
\end{lem}

\proof
Let $\nabla_0$ be a unitary connection of $E_0$ which is
trivial with respect to the trivialization $\tau$ along the boundary
$\partial\Sigma_0$. Over each $(\widehat{D_i}\times\C^n,\Z_{m_i})$,
we define an equivariant connection $\nabla=\beta\psi_i^\ast
\nabla_0+(1-\beta)d$ where $\beta$ is an equivariant cut-off function
equaling one near $\partial\widehat{D_i}$ and $d$ is the trivial
connection with respect to the natural trivialization of
$\widehat{D_i}\times\C^n$. Clearly $\nabla_0,\nabla$ are compatible
on the overlaps so that they define a connection of the orbifold
bundle $E$, which is still denoted by $\nabla$ for simplicity. We
observe that over $\Sigma_0$, $\nabla=\nabla_0$, and with respect to
each local trivialization $(\widehat{D_i}\times\C^n,\Z_{m_i})$, the
curvature form $F(\nabla)$ is given by the diagonal matrix whose entries
are $-d(\beta m_{i,1}\frac{dz}{z}), \cdots, -d(\beta m_{i,n}\frac{dz}{z})$.
Hence
\begin{align*}
c_1(E)(\Sigma) & =  \int_{\Sigma}\frac{\sqrt{-1}}{2\pi}tr F(\nabla)\\
               & =  \int_{\Sigma_0}\frac{\sqrt{-1}}{2\pi}tr F(\nabla_0)
                     + \sum_{i=1}^k\frac{1}{m_i}\int_{\widehat{D_i}}
                       \frac{\sqrt{-1}}{2\pi}tr F(\nabla)\\
               & =  c_1(E_0,\tau)(\Sigma_0,\partial\Sigma_0)+\sum_{i=1}^k
                     (\sum_{j=1}^n\frac{m_{i,j}}{m_i}).\tag*{\qed}
\end{align*}

As an example which is also relevant in the later discussion, we
consider the case where $E=T\Sigma\rightarrow\Sigma$. On each local
uniformizing system $(\widehat{D_i},\Z_{m_i})$, $T\Sigma$ has a natural
trivialization $(\widehat{D_i}\times\C,\Z_{m_i})$ defined by the section
$\frac{\partial}{\partial z}$, where $\Z_{m_i}$ acts by complex
multiplication (ie $m_{i,1}=1$). On the other hand, the trivialization
$\tau$ is defined by $d\phi_i(z\frac{\partial}{\partial z})
=m_i w\frac{\partial}{\partial w}$ along each $\partial D_i$, where
$\phi_i\co \widehat{D_i}\rightarrow D_i$ is the map $w=z^{m_i}$. It is easily
seen that $c_1(T\Sigma_0,\tau)(\Sigma_0,\partial\Sigma_0)=2-2g_{|\Sigma|}-k$
where $g_{|\Sigma|}$ is the genus of the underlying Riemann surface
of $\Sigma$, and $k$ is the number of components in $\partial\Sigma_0$. 
Hence Lemma 3.6 recovers the formula
$$
c_1(T\Sigma)(\Sigma)=2-2g_{|\Sigma|}-\sum_{i=1}^k(1-\frac{1}{m_i}).
$$
Note that the right hand side of the above equation equals $2-2g_\Sigma$
by the definition of the orbifold genus $g_\Sigma$.

\sh{Proof of Theorem 3.1}

We consider first the case where $C$ is a type I $J$--holomorphic curve.
We shall continue to use the notations introduced in the proof of Lemma 3.5.

Let $E\rightarrow\Sigma$ be the pullback of $TX$ by $f$, which is a rank $2$
orbifold complex vector bundle. Over each local uniformizing system
$(\widehat{D_i},\Z_{m_i})$, $E$ has a trivialization $(\widehat{D_i}\times
\C^2,\Z_{m_i})$, where $\{z\}\times\C^2, \forall z\in\widehat{D_i}$, is
identified with $T\widehat{V_i}|_{f_i(z)}$, and $\Z_{m_i}$ acts by
$\mu_{m_i}\cdot (z,w)=(\mu_{m_i}z,\rho_i(\mu_{m_i})(w))$, $\mu_{m_i}=
\exp(\sqrt{-1}\frac{2\pi}{m_i})$. More concretely, we may identify
$\widehat{V_i}$ with $\C^2$ such that the almost complex structure $J$
equals the standard one at the origin $0$, and there are coordinates
$u,v$ such that $\rho_i(\mu_{m_i})$ acts linearly as a diagonal matrix,
say with entries $\mu_{m_i}^{m_{i,1}}$, $\mu_{m_i}^{m_{i,2}}$ where
$0\leq m_{i,1},m_{i,2}<m_i$, and that
$f_i(z)=(z^{l_i},a_iz^{l_i})+O(|z|^{l_i+1})$ for some integer $l_i\geq 1$
and $a_i\in\C$. Observe that if $a_i\neq 0$, then $f_i$ being
$\rho_i$--equivariant implies that $m_{i,1}=m_{i,2}$, so that we may
modify with a linear coordinate change $(u,v)\mapsto (u,v-a_i u)$ such that
$\rho_i(\mu_{m_i})$ is still diagonalized and $f_i(z)=(z^{l_i},0)
+O(|z|^{l_i+1})$. Thus in any event, we have $f_i(z)=(z^{l_i},0)+
O(|z|^{l_i+1})$. Let $E_0=E|_{\Sigma_0}$, and $\tau$ be the canonical
trivialization of $E_0$ along $\partial\Sigma_0$ which is determined by the
equivariant sections $(z^{m_{i,1}},0)$ and $(0,z^{m_{i,2}})$ of
$\widehat{D_i}\times\C^2\rightarrow\widehat{D_i}$ along each $\partial
D_i$. Recall that $c=-c_1(TX)$. Hence by Lemma 3.6,
$$
c(C)=-c_1(E_0,\tau)(\Sigma_0,\partial\Sigma_0)-\sum_{i=1}^k
\frac{m_{i,1}+m_{i,2}}{m_i}.
$$

Observe that $f_\epsilon^\ast TX^0=E_0$ along $\partial\Sigma_0\subset
\Sigma_0^\prime$. Hence the canonical trivialization $\tau$ of $E_0$ along
$\partial\Sigma_0$ gives rise to a trivialization of $f_\epsilon^\ast TX^0$
along $\partial\Sigma_0^\prime$, which is also denoted by $\tau$ for
simplicity. Furthermore, note that $c_1(E_0,\tau)(\Sigma_0,\partial\Sigma_0)
=c_1(f^\ast_\epsilon TX^0,\tau)(\Sigma_0^\prime,\partial\Sigma_0^\prime)$.
On the other hand, let $\tau_h$ be the trivialization of $T\Sigma_0^\prime$
along the
boundary $\partial\Sigma_0^\prime$ given by the section $w\frac{\partial}
{\partial w}$ (here $w$ is the holomorphic coordinate of each
$D_i$). Then $\tau,\tau_h$ determine a unique trivialization $\tau_v$
of $\nu_\epsilon$ along $\partial\Sigma_0^\prime$ such that
$$
c_1(f_\epsilon^\ast TX^0,\tau)=c_1(T\Sigma_0^\prime,\tau_h)+
c_1(\nu_\epsilon,\tau_v).
$$

There are canonical bundle morphisms
$\nu_{i,\epsilon}|_{\partial \widehat{D_i}}\rightarrow
\nu_\epsilon|_{\partial D_i}$ induced by $\phi_i\co \widehat{D_i}
\rightarrow D_i$ where $\phi_i(z)=z^{m_i}$. Through these bundle morphisms,
the trivialization $\tau_v$ gives rise to a trivialization
$\tau_{i,v}$ of $\nu_{i,\epsilon}$ along $\partial \widehat{D_i}$.
In order to determine $\tau_{i,v}$, we recall that $f_i(z)=(z^{l_i},0)
+O(|z|^{l_i+1})$ and $f_{i,\epsilon}=f_i$ in $\widehat{D_i}\setminus
\widehat{D_i^\prime}$. If we let $\tau_{i,h}$ be the trivialization of
$T\widehat{D_i}$ along $\partial\widehat{D_i}$ (as a sub-bundle of
$f_{i,\epsilon}^\ast T\widehat{V_i}$) which is induced by the
trivialization $\tau_h$ of $T\Sigma_0^\prime$ along $\partial\Sigma_0^\prime$
through $\phi_i$, then $\tau_{i,h}$ is given by
the section $(l_i z^{l_i},0)$ up to homotopy. Hence $\tau_{i,v}$
is given by the section $(0,z^{-l_i+m_{i,1}+m_{i,2}})$ up to
homotopy, since $\tau$ is given by the sections $(z^{m_{i,1}},0)$
and $(0,z^{m_{i,2}})$.

We push $f_\epsilon$ off near $\partial\Sigma_0^\prime$ along the
direction given by the trivialization $\tau_v$ of the normal
bundle $\nu_\epsilon$ (note that $f_\epsilon$ is embedded near
$\partial\Sigma_0^\prime$). Call the resulting map
$f_\epsilon^\prime$. Correspondingly, each $f_{i,\epsilon}$ is
pushed off near $\partial\widehat{D_i}$ to a $f_{i,\epsilon}^\prime$
along the direction given by the trivialization $\tau_{i,v}$ of the
normal bundle $\nu_{i,\epsilon}$. As in the proof of Lemma 3.5, we can
similarly construct a closed 2--form $\eta_C^\prime$ using $f_\epsilon^\prime,
f_{i,\epsilon}^\prime$ instead of $f_\epsilon, f_{i,\epsilon}$, which is
also Poincar\'{e} dual to the differentiable cycle $f(\Sigma)$. Furthermore,
$$
C\cdot C = \int_X \eta_C^\prime\wedge\eta_C =
\int_{\Sigma_0}(f_\epsilon^\prime)^\ast\eta_C+
               \sum_{i=1}^k\frac{1}{m_i}\int_{\widehat{D_i}}
               (f_{i,\epsilon}^\prime)^\ast\eta_C.
$$
By way of construction,
$$
\int_{\Sigma_0}(f_\epsilon^\prime)^\ast\eta_C=
c_1(\nu_\epsilon,\tau_v)(\Sigma_0,\partial\Sigma_0)+
\sum_{\{[z,z^\prime]|z\neq z^\prime,
f_\epsilon(z)=f_\epsilon(z^\prime)\}}2t_{[z,z^\prime]},
$$
where $[z,z^\prime]$ denotes the unordered pair of $z,z^\prime$, and
$t_{[z,z^\prime]}$ is the order of tangency of the
intersection $f_\epsilon(z)=f_\epsilon(z^\prime)$. It is easily
seen that the second term in the above equation is equal to
$$
\sum_{\{[z,z^\prime]|z\neq z^\prime,f(z)=f(z^\prime)\}}2k_{[z,z^\prime]} +
\sum_{\{z|df(z)=0\}}2k_z, \mbox{  where  } z,z^\prime\in\Sigma_0.
$$

To evaluate $\int_{\widehat{D_i}}(f_{i,\epsilon}^\prime)^\ast\eta_C$,
$i=1,2,\cdots,k$, let $I_i$ be the set labeling $\Lambda(C)_{z_i}$, ie
$\Lambda(C)_{z_i}=\{C_{i,\alpha}\mid\alpha\in I_i\}$, and let 
$C_i\in\Lambda(C)_{z_i}$ be the element defined by $f_i$. Then
\begin{align*}
\int_{\widehat{D_i}}(f_{i,\epsilon}^\prime)^\ast\eta_C & = 
\int_{\widehat{D_i}}(f_{i,\epsilon}^\prime)^\ast(\sum_{\{j|
f(z_i)=f(z_j)\}}\frac{1}{m_j}\sum_{g\in
G_i}g^\ast\overline{\Theta_{j,\epsilon}})\\
& =  \frac{1}{m_i}\sum_{g\in G_i}\int_{\widehat{D_i}}
(f_{i,\epsilon}^\prime)^\ast
(g^\ast\overline{\Theta_{i,\epsilon}})\\
&\hspace{1.5in} + \sum_{\{j\neq
i|f(z_j)=f(z_i)\}}\frac{1}{m_j}\sum_{g\in G_i}\int_{\widehat{D_i}}
(f_{i,\epsilon}^\prime)^\ast(g^\ast\overline{\Theta_{j,\epsilon}})\\
& =  c_1(\nu_{i,\epsilon},\tau_{i,v})(\widehat{D_i},\partial\widehat{D_i})+
C_{i}\cdot C_{i}+\sum_{\alpha\in I_i}C_{i}\cdot
C_{i,\alpha}\\
&\hspace{1.5in}+ \sum_{\{j\neq i|f(z_i)=f(z_j)\}}\sum_{\beta\in I_j}
C_{i}\cdot C_{j,\beta}\\
& =  c_1(\nu_{i,\epsilon},\tau_{i,v})(\widehat{D_i},\partial\widehat{D_i})+
\frac{m_i}{|G_i|}(\sum_{\alpha\in I_i}C_{i,\alpha}\cdot C_{i,\alpha}+
\sum_{\alpha,\beta\in I_i}C_{i,\alpha}\cdot C_{i,\beta})\\
&\hspace{1.5in}   + \frac{m_i}{|G_i|}\sum_{\{j\neq i|f(z_i)=f(z_j)\}}
\sum_{\alpha\in I_i,\beta\in I_j}C_{i,\alpha}\cdot C_{j,\beta}.
\end{align*}
In order to evaluate $c_1(\nu_{i,\epsilon},\tau_{i,v})
(\widehat{D_i},\partial\widehat{D_i})$, we observe that $f_{i,\epsilon}$
is an immersion and equals $(z^{l_i},0)+O(|z|^{l_i+1})$ near $\partial
\widehat{D_i}$. Let $\tau_{i,v}^\prime$ be the trivialization of
$\nu_{i,\epsilon}$ along $\partial\widehat{D_i}$ which can be extended
over the entire $\widehat{D_i}$. Then $\tau_{i,v}^\prime$ is given by
the section $(0,z^{-l_i+1})$ up to homotopy. But $\tau_{i,v}$ is given
by the section $(0,z^{-l_i+m_{i,1}+m_{i,2}})$ up to homotopy. Hence
$$c_1(\nu_{i,\epsilon},\tau_{i,v})(\widehat{D_i},
\partial\widehat{D_i})=m_{i,1}+m_{i,2}-1.$$ Putting things altogether, we have
\begin{eqnarray*}
C\cdot C + c(C) & = & c(C)+ c_1(\nu_\epsilon,\tau_v)(\Sigma_0,
\partial\Sigma_0)+\sum_{i=1}^k\frac{m_{i,1}+m_{i,2}-1}{m_i} \\
&  & + \sum_{\{[z,z^\prime]|z\neq z^\prime,
f(z)=f(z^\prime)\}}2k_{[z,z^\prime]} +
\sum_{z\in\Sigma}2k_z\\
& = & -c_1(T\Sigma_0,\tau_h)(\Sigma_0,\partial\Sigma_0)-\sum_{i=1}^k\frac{1}
{m_i}\\
&  & + \sum_{\{[z,z^\prime]|z\neq z^\prime,
f(z)=f(z^\prime)\}}2k_{[z,z^\prime]} +
\sum_{z\in\Sigma}2k_z\\
& = & 2g_{|\Sigma|}-2+k-\sum_{i=1}^k\frac{1}{m_i}\\
&  & + \sum_{\{[z,z^\prime]|z\neq z^\prime,
f(z)=f(z^\prime)\}}2k_{[z,z^\prime]} +
\sum_{z\in\Sigma}2k_z,
\end{eqnarray*}
from which the adjunction formula for the case where $C$ is of type I
follows easily.

The case where $C$ is of type II is actually much simpler. It follows
by directly evaluating the last integral in
$$
C\cdot C=\int_X\eta_C\wedge\eta_C=\int_\Sigma f^\ast\eta_C,
$$
and then appealing to $c_1(TX)(\Sigma)=c_1(\nu)(\Sigma)+c_1(T\Sigma)(\Sigma)$
and $m_C=m_\Sigma$.
\endproof

\sh{Proof of Theorem 3.2}

For simplicity, we shall only consider the case where $C,C^\prime$
are of type I. The discussion for the rest of the cases is similar, and
we shall leave the details to the reader.

Let $\eta_C$, $\eta_{C^\prime}$ be the closed 2--forms in Lemma 3.5 which
are Poicar\'{e} dual to the differentiable cycles $f(\Sigma)$,
$f^\prime(\Sigma^\prime)$ respectively. Then
\begin{eqnarray*}
C\cdot C^\prime & = & \int_X\eta_C\wedge\eta_{C^\prime}\\
                & = &
                \int_{\Sigma_0}f^\ast_\epsilon\eta_{C^\prime}+
                \sum_{i=1}^k\frac{1}{m_i}\int_{\widehat{D_i}}
                f^\ast_{i,\epsilon}\eta_{C^\prime}.
\end{eqnarray*}
Now observe that the subset
$\{(z,z^\prime)\mid f(z)=f^\prime(z^\prime)\}\subset\Sigma\times\Sigma^\prime$
is finite. Hence we may arrange in the construction of $\eta_C$ and
$\eta_{C^\prime}$ such that for sufficiently small $\epsilon>0$,
$\int_{\Sigma_0}f^\ast_\epsilon\eta_{C^\prime}$ equals $\sum k_{(z,z^\prime)}$
where $(z,z^\prime)$ is running over the set of pairs with
$f(z)=f^\prime(z^\prime)$ being a regular point of $X$, and
$\sum_{i=1}^k\frac{1}{m_i}\int_{\widehat{D_i}}f^\ast_{i,\epsilon}
\eta_{C^\prime}$ equals $\sum k_{(z,z^\prime)}$ where $(z,z^\prime)$ is
running over the set of pairs with $f(z)=f^\prime(z^\prime)$
being an orbifold point of $X$. Hence the theorem.
\endproof

\section{Proof of main results}

We begin by setting the stage. Let $p,q$ be relatively prime integers
with $0<q<p$. We denote by $C_{(p,q)}$ the symplectic cone over $L(p,q)$,
which is the symplectic orbifold $(\C^2,\omega_0)/\Z_p$ where
$\omega_0=\frac{\sqrt{-1}}{2}\sum_{i=1}^2 dz_i\wedge d\bar{z_i}$ and
$\Z_p$ acts by $\mu_p\cdot(z_1,z_2)=(\mu_p z_1,\mu_p^q z_2)$. Let $d$ be
the descendant of the function $\frac{1}{2}(|z_1|^2+|z_2|^2)$ on $\C^2$ to
$C_{(p,q)}$. Then for any $r>0$, $C_{(p,q)}(r)\equiv d^{-1}([0,r])
\subset C_{(p,q)}$ is a suborbifold of contact boundary $(L(p,q),\xi_0)$.

Next we follow the discussion in \cite{Le} to embed each $C_{(p,q)}(r)$
into an appropriate closed symplectic 4--orbifold. To this end, consider
the Hamiltonian circle action on $(\C^2,\omega_0)$
$$
s\cdot (z_1,z_2)=(sz_1,s^{p+q} z_2), \forall s\in\s^1\equiv \{z\in\C\mid
|z|=1\},
$$
with the Hamiltonian function given by $\mu(z_1,z_2)
=\frac{1}{2}(|z_1|^2+(p+q)|z_2|^2)$. It is easily seen that the
$\Z_p$--action on $\C^2$ is the action induced from the circle action by
$\Z_p\subset\s^1$, thus there is a corresponding Hamiltonian circle
action on $\C^2/\Z_p=C_{(p,q)}$ with the Hamiltonian function given by
$\mu^\prime\equiv\frac{1}{p}\mu$. According to \cite{Le}, for any $R>0$,
there is a symplectic 4--orbifold, denoted by $X_{(p,q)}(R)$, which is
obtained from $(\mu^\prime)^{-1}([0,R])$ by collapsing each orbit of
the circle action on $(\mu^\prime)^{-1}(R)$ to a point. It is clear that
for any $R>\frac{1}{p}(p+q)r$, $C_{(p,q)}(r)$ is a suborbifold of
$X_{(p,q)}(R)$ of contact boundary $(L(p,q),\xi_0)$. Furthermore, there is
a distinguished 2--dimensional symplectic suborbifold
$C_0\equiv (\mu^\prime)^{-1}(R)/\s^1\subset X_{(p,q)}(R)$, whose normal
bundle has Euler number $\frac{p}{p+q}$, and whose orbifold genus is
$\frac{1}{2}-\frac{1}{2(p+q)}$, cf Section 3.

Now let $(W,\omega)$ be a symplectic cobordism from
$(L(p^\prime,q^\prime),\xi_0^\prime)$ to $(L(p,q),\xi_0)$. By adding
appropriate ``symplectic collars'' to the two ends of $W$, which does not
change the diffeomorphism class of $W$, we may assume 
without loss of generality that a neighborhood of $L(p^\prime,q^\prime)$ 
in $W$ is identified with a neighborhood of
$\partial C_{(p^\prime,q^\prime)}(r^\prime)$ in $C_{(p^\prime,q^\prime)}
\setminus \mbox{int}(C_{(p^\prime,q^\prime)}(r^\prime))$ for some
$r^\prime>0$, and a neighborhood of $L(p,q)$ in $W$ is 
identified with a neighborhood of $\partial C_{(p,q)}(r)$ in $C_{(p,q)}(r)$
for some $r>0$. Consequently, we can close up $W$ by gluing
$X_{(p,q)}(R)\setminus C_{(p,q)}(r)$ and $C_{(p^\prime,q^\prime)}(r^\prime)$
onto the corresponding ends of $W$ for some fixed $R>\frac{1}{p}(p+q)r$.
We denote by $(X,\omega)$ the resulting symplectic 4--orbifold. Note that
there is a distinguished 2--dimensional symplectic suborbifold $C_0\subset X$
inherited from $C_0\subset X_{(p,q)}(R)$. 

With the preceding understood, the strategy for proving Theorem 1.1 is
to construct a diffeomorphism of orbifold pairs from $(X,C_0)$ 
to $(X_{(p,q)}(R),C_0)$. 

First of all, some preliminary information about $(X,C_0,\omega)$.
The orbifold $X$ has two singular points, one of them, denoted by
$x^\prime$, is inherited from $C_{(p^\prime,q^\prime)}(r^\prime)$
and has type $(p^\prime,q^\prime)$, and the other, denoted by $x$,
is inherited from $X_{(p,q)}(R)\setminus C_{(p,q)}(r)$ and has
type $(p+q,p)$. Here a singular point has type $(a,b)$ if the isotropy
group is $\Z_a$ with action on a local uniformizing system given by
$\mu_a\cdot (z_1,z_2)=(\mu_a z_1,\mu_a^b z_2)$. The suborbifold $C_0$
has only one orbifold point, the point $x$ with order $p+q$, and is
given locally by $z_2=0$ on the local uniformizing system. We fix
an $\omega$--compatible almost complex structure $J$ on $X$, such that
the suborbifold $C_0$ is $J$--holomorphic. For convenience, we assume
that $J$ is integrable near $x,x^\prime$. (This is possible because
of the equivariant Darboux' theorem.) By the discussion in Section 3,
we see that $C_0\cdot C_0=e(\nu)(C_0)=\frac{p}{p+q}$, where $e(\nu)$ is the
Euler class of the normal bundle $\nu$ of $C_0$ in $X$. On the other hand,
by the adjunction formula in Theorem 3.1, we have
$$
c_1(K_X)(C_0)=2(\frac{1}{2}-\frac{1}{2(p+q)})-2-C_0\cdot C_0
=-\frac{2p+q+1}{p+q}
$$
for the canonical bundle $K_X$ of the almost complex 4--orbifold $(X,J)$.

Next we digress on the Fredholm theory for pseudoholomorphic curves in a
symplectic 4--orbifold $(X,\omega)$. To this end, for any given orbifold
Riemann surface $\Sigma$, we fix a sufficiently large positive integer $k$,
and consider $[\Sigma;X]$, the space of $C^k$ maps from $\Sigma$ into $X$.
It is shown in \cite{C1} (Part I, Theorem 1.4) that $[\Sigma;X]$ is a smooth
Banach orbifold (Hausdorff and second countable). Moreover, a map
$f\in [\Sigma;X]$ is a smooth point in the Banach orbifold if $\mbox{Im }f$
contains a regular point of $X$. Thus for the purpose here we may assume for
simplicity that $\Sigma$ is reduced and $[\Sigma;X]$ is a Banach manifold.
The tangent space $T_f$ at $f\in [\Sigma;X]$ is the space of $C^k$ sections
of $f^\ast(TX)$, the pullback bundle of $TX$ via $f$.

For any $f\in [\Sigma;X]$, let $\E_f$ be the subspace of the space of
$C^{k-1}$ sections of the orbifold vector bundle 
$\mbox{Hom}(T\Sigma,f^\ast(TX))\rightarrow\Sigma$, which consists of
sections $s$ satisfying $s\circ j=-J\circ s$ for a fixed choice of 
$\omega$--compatible almost complex structure $J$ on $X$ and the complex 
structure $j$ on $\Sigma$. Then there is a Banach bundle $\E$ over 
$[\Sigma;X]$ whose fiber at $f$ is $\E_f$. Consider the smooth section 
$\underline{L}\co [\Sigma;X]\rightarrow\E$ defined by
$$
\underline{L}(f)\equiv df+J\circ df\circ j.
$$
The zero loci $\underline{L}^{-1}(0)$ is the space of $J$--holomorphic
maps from $\Sigma$ into $X$. By elliptic regularity, each map in
$\underline{L}^{-1}(0)$ is a $C^\infty$ map. Moreover, $\underline{L}$
is a Fredholm section, and its linearization $D\underline{L}$ at each
$f\in \underline{L}^{-1}(0)$ is given by a formula
$$
D\underline{L}_f(u)=L_{f}(u),\hspace{3mm} u\in T_f,
$$
where $L_f\co T_f\rightarrow\E_{f}$ is an elliptic linear differential
operator of Cauchy--Riemann type, whose coefficients are smooth functions
on $\Sigma$ which depend on $f$ smoothly. The following facts are crucial 
for the consideration of surjectivity of $D\underline{L}$. 
\begin{itemize}
\item When $J$ is integrable in a neighborhood of $\mbox{Im }f$ and
$f$ is $J$--holomorphic, $D\underline{L}_f=L_f$ is the usual
$\bar{\partial}$--operator for the orbifold holomorphic vector bundle
$f^\ast (TX)$ over $\Sigma$.
\item When $f$ is a multiplicity-one parametrization of a $J$--holomorphic
suborbifold $C$, the linearization $D\underline{L}_f=L_f$ is surjective
when $c_1(TC)(C)>0$ and $c_1(K_X)(C)<0$. This is the orbifold analog of
the regularity criterion discussed in Lemma 3.3.3 of \cite{McDS}.
\end{itemize}

The index of $D\underline{L}_f=L_{f}$ can be computed using the index
formula of Kawasaki \cite{Ka} for elliptic operators on orbifolds, cf
Lemma 3.2.4 in \cite{CR}.\par To state the formula, let $z_1,z_2,\cdots,z_l$
be the set of orbifold points of $\Sigma$ with orders $m_1,m_2,\cdots,
m_l$ respectively. Moreover, suppose at each $z_i$, a local representative
of $f$ is given by $(f_i,\rho_i)\co (\widehat{D_i},\Z_{m_i})\rightarrow
(\widehat{V_i},G_i)$ where $\rho_i(\mu_{m_i})$ acts on $\widehat{V_i}$
by $\rho_i(\mu_{m_i})\cdot (w_1,w_2)=(\mu_{m_i}^{m_{i,1}}w_1,
\mu_{m_i}^{m_{i,2}}w_2)$, $0\leq m_{i,1}, m_{i,2}<m_i$. With this understood,
$\mbox{Index }D\underline{L}_f=2d$ where $d\in\Z$ is given by
$$
d=c_1(TX)\cdot [f(\Sigma)]+2-2g_{|\Sigma|}
-\sum_{i=1}^l\frac{m_{i,1}+m_{i,2}}{m_i}.
$$
(Here $g_{|\Sigma|}$ is the genus of the underlying Riemann surface.)
End of digression.

Now let $\Sigma$ be the orbifold Riemann sphere with one orbifold point
$z_\infty\equiv\infty$ of order $p+q$. Observe that as a complex
analytic space, $\Sigma$ is biholomorphic to the underlying
Riemann sphere $|\Sigma|$, hence it has a unique complex
structure. Moreover, the group of automorphisms $G$ can be
naturally identified with the subgroup of the automorphism group
of $|\Sigma|$ which fixes the point $\infty$. Note that $|\Sigma|\setminus
\{\infty\}=\C$, so that $G$ can be identified with the group $\{(a,b)
\in C^\ast\times\C\mid z\mapsto az+b\}$ of linear translations on $\C$.

We shall consider the moduli space $\widetilde{\M}$ of $J$--holomorphic maps
$f\co \Sigma\rightarrow X$ which obey
\begin{itemize}
\item $[f(\Sigma)]=[C_0]$ in $H_2(X;\Q)$, 
\item $f(z_\infty)=x$, and in a local representative 
$(f_\infty,\rho_\infty)$ of $f$ at $z_\infty$, $\rho_\infty(\mu_{(p+q)})$
$=\mu_{(p+q)}$, which acts by $(z_1,z_2)\mapsto (\mu_{(p+q)} z_1,
\mu_{(p+q)}^p z_2)$. (Here $z_1,z_2$ are holomorphic coordinates 
on a local uniformizing system at $x$ in which $C_0$ is locally given 
by $z_2=0$.)
\end{itemize}
We set $\M=\widetilde{\M}/G$ for the corresponding moduli space of
unparametrized $J$--holomorphic maps, where $G$ acts on $\widetilde{\M}$
by reparametrization.

With the preceding understood, consider the following:

\begin{lem}
Suppose $W$ is a {\em(}symplectic{\em)} homology cobordism. 
{\em(}Note that in particular, 
$p=p^\prime$ and $H_2(X;\Q)=\Q\cdot [C_0]$.{\em)} Then 

{\em (1)}\qua Each member of $\widetilde{\M}$ is either an orbifold embedding
onto a suborbifold in $X$, or is a multiply covered map with multiplicity 
$p$ onto a suborbifold containing both $x,x^\prime$. Moreover, in the 
latter case, either $q^\prime=q$ or $q^\prime q\equiv 1\pmod{p}$ must be
satisfied, and there is at most one such a member of $\widetilde{\M}$ up
to reparametrization by elements of $G$.

{\em (2)}\qua One may alter $J$ appropriately such that $C_0$ is still
$J$--holomorphic, and $\widetilde{\M}$ is a smooth manifold of
dimension $6$. Furthermore, $\M$ is a compact, closed, $2$--dimensional
smooth orbifold (possibly disconnected) with at most one orbifold point
of order $p$, and the action of $G$ on $\widetilde{\M}$ defines a smooth
orbifold principal $G$--bundle $\widetilde{\M}\rightarrow\M$.
\end{lem}

Before proving Lemma 4.1, let us observe the following:

\begin{lem}
Let $C$ be any $J$--holomorphic curve in $X$ such that
\begin{itemize}
\item $C$ contains both singular points,
\item $[C]=r[C_0]$ for some $r\in (0,1]\cap\Q$. 
\end{itemize}
Then $C$ is a suborbifold and $[C]=\frac{1}{p}[C_0]$. Moreover,
there is at most one such $J$--holomorphic curves in $X$. 
\end{lem}

\proof
First of all, we claim $r\geq\frac{1}{p}$. To see this, note that 
$C\neq C_0$ because $C$ contains both singular points. By the intersection
formula (cf Theorem 3.2), 
$$
r\cdot\frac{p}{p+q}=C\cdot C_0\geq\frac{1}{p+q},
$$ 
which verifies the claim. 

Now let $f\co \Sigma\rightarrow X$ be a multiplicity-one parametrization 
of $C$, and $z_0,z_0^\prime\in \Sigma$ be any points such that $f(z_0)=x$,
$f(z_0^\prime)=x^\prime$. Let $m_0,m_0^\prime$ be the order of 
$z_0,z_0^\prime$ respectively. Then observe that if $m_0<p+q$ (resp.
$m_0^\prime<p$), the contribution $k_{z_0}$ (resp. $k_{z_0^\prime}$)
on the right hand side of the adjunction formula for $C$ (cf Theorem 3.1)
is no less than $\frac{1}{2m_0}$ (resp. $\frac{1}{2m_0^\prime}$).
(Here is the calculation for the case of $m_0$: $k_{z_0}\geq \frac{1}{2(p+q)}
\cdot [\frac{p+q}{m_0}(\frac{p+q}{m_0}-1)]\geq \frac{1}{2m_0}$ if $m_0<p+q$.)
It follows easily that the right hand side of the adjunction formula for $C$
is no less than 
$$
\frac{1}{2}(1-\frac{1}{p+q})+\frac{1}{2}(1-\frac{1}{p}),
$$
which has an equality only if $m_0=p+q$ and $m_0^\prime=p$. 

On the other hand, the left hand side of the adjunction formula for $C$,
the virtual genus $g(C)$, equals 
$$
\frac{1}{2}(\frac{p}{p+q}\cdot r^2-\frac{2p+q+1}{p+q}\cdot r)+1.
$$
As a function of $r$, it is decreasing over $(0,1]$, hence the maximum
of $g(C)$ is attained at $r=\frac{1}{p}$, and it equals
$$
\frac{1}{2}(\frac{p}{p+q}\cdot (\frac{1}{p})^2-\frac{2p+q+1}{p+q}\cdot 
\frac{1}{p})+1=\frac{1}{2}(1-\frac{1}{p+q})+\frac{1}{2}(1-\frac{1}{p}).
$$
By the adjunction formula, $C$ is a suborbifold and $[C]=\frac{1}{p}[C_0]$.

To see that there is at most one such $J$--holomorphic curves, note 
that if there were two distinct such curves, the algebraic intersection
number, which is $\frac{p}{p^2(p+q)}$, would be at least
$\frac{1}{p+q}+\frac{1}{p}$ by the intersection formula.
A contradiction.

\sh{Proof of Lemma 4.1}

(1)\qua By the adjunction formula, each multiplicity-one
member $f\in\widetilde{\M}$ must be an orbifold embedding
onto a suborbifold. Now suppose 
$f\in\widetilde{\M}$ is multiply covered with multiplicity $m>1$. Let
$C$ be the corresponding $J$--holomorphic curve. Then $[C]=\frac{1}{m}[C_0]
<[C_0]$, which implies that $C$ also contains the other singular point 
$x^\prime$. This is because by the assumption, $W$ is a homology
cobordism, so that $H_2(X\setminus\{x^\prime\};\Z)$ is generated by 
the class of $C_0$, and hence $C$ can not be contained entirely in
$X\setminus\{x^\prime\}$. By Lemma 4.2, $f$ has multiplicity $p$, and 
$C$ is a suborbifold, which is unique in such kind. 

To complete the proof of (1), it remains to show that either $q^\prime=q$ 
or $q^\prime q\equiv 1\pmod{p}$ if there is indeed such a curve $C$. 

To this end, let $f\co \Sigma\rightarrow X$ be any multiplicity-one 
parametrization of $C$, and $z_0,z_0^\prime\in \Sigma$ be the points 
such that $f(z_0)=x$, $f(z_0^\prime)=x^\prime$. Since $C\neq C_0$ and
$C\cdot C_0=\frac{1}{p}\cdot\frac{p}{p+q}=\frac{1}{p+q}$, it follows
easily that the local representative of $f$ at $z_0$ must be in the form 
$((u(z),z),\rho_0)$ for some holomorphic function $u$ and the 
isomorphism $\rho_0$ where $\rho_0(\mu_{(p+q)})=\mu_{(p+q)}^l$ with
$pl\equiv 1\pmod{p+q}$. On the other hand, the local representative of 
$f$ at $z_0^\prime$ could either be $((w(z),z),\rho_0^\prime)$, where
$\rho_0^\prime(\mu_{p})=\mu_{p}^{l^\prime}$ with $l^\prime q^\prime
\equiv 1 \pmod{p}$, or $((z,w(z)),\rho_0^\prime)$ with
$\rho_0^\prime(\mu_{p})=\mu_{p}$. Assuming the former case, we have,
by the index formula for $D\underline{L}_f$,
$$
\frac{2p+q+1}{p(p+q)}+2-\frac{l+1}{p+q}-\frac{l^\prime+1}{p}\in\Z,
$$
which implies that $r(p+q)-ql^\prime\equiv 0\pmod{p}$ with $r$ given
by the equation $1-lp=r(p+q)$. It is easily seen that in this case,
$ql^\prime\equiv qr\equiv 1\pmod{p}$, and hence $q^\prime=q$ because
$l^\prime q^\prime\equiv 1 \pmod{p}$. Similarly,
the latter case implies $q^\prime q\equiv 1 \pmod{p}$.

(2)\qua For the smoothness of $\widetilde{\M}$, we need to show that
for any $f\in\widetilde{\M}$, the linearization $D\underline{L}_f$
is surjective. The dimension of $\widetilde{\M}$ is the index of
$D\underline{L}_f$, $f\in\widetilde{\M}$, which is easily seen to
be $6$ by the index formula for $D\underline{L}_f$.

By the regularity criterion we mentioned earlier, $\widetilde{\M}$ is 
smooth at each
$f$ which is not multiply covered, because for any such an $f$, $C\equiv
\mbox{Im} f$ is a suborbifold satisfing $c_1(TC)(C)=2-(1-\frac{1}{p+q})>0$
and $c_1(K_X)(C)=-\frac{2p+q+1}{p+q}<0$. Suppose there is a multiply
covered member (which is the only one up to reparametrization by (1)),
and let $C_0^\prime$ be the corresponding $J$--holomorphic curve. We consider
the weighted projective space $\P(1,p,p+q)$, which is the quotient
of $\s^5$ under the $\s^1$--action
$$
s\cdot (z_1,z_2,z_3)=(s z_1,s^p z_2,s^{p+q} z_3), \forall s\in\s^1.
$$
It is easily seen that a regular neighborhood of $C_0^\prime$ in $X$
is diffeomorphic to a regular neighborhood of $\P(p,p+q)$ in $\P(1,p,p+q)$,
where $\P(p,p+q)$ is defined by $z_1=0$. According to \cite{BGN},
$\P(1,p,p+q)$ has an orbifold K\"{a}hler metric of positive Ricci curvature.
By the orbifold version of symplectic neighborhood theorem, we can
alter the almost complex structure $J$ in a regular neighborhood of
$C^\prime_0$ such that $\omega(\cdot,J(\cdot))$ is K\"{a}hler of positive
Ricci curvature. (Note that we can arrange so that $C_0$ is still
$J$--holomorphic, and $J$ is integrable near singular points $x,x^\prime$.)
With this understood, for any $f\in\widetilde{\M}$ parametrizing $C_0^\prime$,
$D\underline{L}_f$ is the usual $\bar{\partial}$--operator for the
orbifold holomorphic vector bundle $f^\ast (TX)$ over $\Sigma$. In this
case, the surjectivity of $D\underline{L}_f$ follows from the orbifold
version of a Bochner type vanishing theorem for negative holomorphic
vector bundles (cf \cite{Ko}). Thus in any event, by altering $J$ if 
necessary, we can arrange so that $\widetilde{\M}$ is a smooth manifold.

The action of $G$ on $\widetilde{\M}$ is smooth (see the general discussion
at the end of \S 3.3 of Part I of \cite{C1}), and is free at each
$f\in\widetilde{\M}$ which is not multiply covered. At a multiply covered
$f\in\widetilde{\M}$, the isotropy subgroup is the cyclic subgroup
$\{(\mu_p^l,0)\mid l=0,\cdots,p-1\}\subset G$ of order $p$ up to
conjugation. (Note that $p$ equals the multiplicity of the covering.) Thus
$\widetilde{\M}\rightarrow\widetilde{\M}/G=\M$ is a smooth orbifold
principal $G$--bundle over a smooth $2$--dimensional orbifold with at
most one orbifold point of order $p$.

It remains to show that $\M$ is compact. First of all, by the orbifold 
version of the Gromov's compactness theorem (cf \cite{Gr, PW, Ye}) 
which was proved in \cite{CR}, any sequence of maps $f_n\in\widetilde{\M}$ 
has a subsequence which converges to a cusp-curve after
suitable reparametrization. More concretely, after reparametrization
if necessary, there is a subsequence of $f_n$, which is still denoted by
$f_n$ for simplicity, and there are at most finitely many simple closed
loops $\gamma_1,\cdots,\gamma_l\subset\Sigma$ containing no orbifold points,
and a nodal orbifold Riemann surface $\Sigma^\prime=\cup_\omega
\Sigma_\omega$ obtained by collapsing $\gamma_1,\cdots,\gamma_l$, and
a $J$--holomorphic map $f\co \Sigma^\prime\rightarrow X$, 
such that (1) $f_n$ converges in $C^\infty$ to $f$ on any given compact
subset in the complement of $\gamma_1,\cdots,\gamma_l$, 
(2) $[f_n(\Sigma)]=[f(\Sigma^\prime)]\in H_2(X;\Q)$, and (3) $f\in
\widetilde{\M}$ and $f_n$ converges to $f$ in $C^\infty$ if there is 
only one component of $\Sigma^\prime=\cup_\omega\Sigma_\omega$ over 
which $f$ is nonconstant.

Hence the space $\M$ is compact if there is only one component of
$\Sigma^\prime=\cup_\omega\Sigma_\omega$ over which $f$ is nonconstant.
Suppose this is not true. Then there is a nonconstant component
$f_\omega\equiv f|_{\Sigma_\omega}\co \Sigma_\omega\rightarrow X$, where
$\Sigma_\omega$ is obtained by collapsing a simple closed
loop $\gamma\in\{\gamma_1,\cdots,\gamma_l\}$ which bounds a disc
$D\subset\Sigma$, such that $z_\infty\in\Sigma\setminus D$ and $f_n$
converges to $f_\omega$ in $C^\infty$ on any compact subset of the interior
of $D$. Set $C_\omega\equiv \mbox{Im } f_\omega$. Since we assume that
there are more than one nonconstant components, $[C_\omega]
\leq [f_\omega(\Sigma_\omega)]<[C_0]$ must hold. (Note that 
$H_2(X;\Q)=\Q\cdot [C_0]$.) By the assumption 
that $W$ is a homology cobordism, $C_\omega$
must contain the singular point $x^\prime$ as we argued earlier. 
We claim that $C_\omega$ must also contain the other singular point.
Suppose not, then $C_\omega\neq C_0$, and $C_\omega$ must intersect
with $C_0$ at a smooth point, because $C_\omega\cdot C\neq 0$. Then
by the intersection formula, $C_\omega\cdot C_0\geq 1$, which implies
that $[C_\omega]=r[C_0]$ for some $r\geq 1+\frac{q}{p}$. A contradiction
to $[C_\omega]<[C_0]$. Now by Lemma 4.2, $C_\omega$ is a suborbifold and
$[C_\omega]=\frac{1}{p}[C_0]$. 

On the other hand, observe that there is a regular point $z_0\in\Sigma_\omega$
such that either $f_\omega(z_0)=x$ or $f_\omega(z_0)=x^\prime$. Let 
$m_\omega\geq 1$ be the multiplicity of $f_\omega$, and let $D_0$
be a sufficiently small disc neighborhood of $z_0$ in $\Sigma_\omega$.
Then it is easily seen that $m_\omega$ is no less than the degree of 
the covering map $f_\omega|_{\partial D_0}$ onto the link of $f_\omega(z_0)$ 
in $C_\omega$, which is no less than $p+q$ or $p$, depending on whether 
$f_\omega(z_0)=x$ or $f_\omega(z_0)=x^\prime$. In any event, 
$m_\omega\geq p$. But this contradicts $[C_\omega]=\frac{1}{p}[C_0]$ as 
$[C_\omega]=\frac{1}{m_\omega}[f_\omega(\Sigma_\omega)]<\frac{1}{p}[C_0]$, 
because $[f_\omega(\Sigma_\omega)]<[C_0]$. 

Hence there is only one nonconstant component, and therefore $\M$ is
compact. 
\endproof

Let $H=\C^\ast$ be the subgroup of $G=\{(a,b)\in\C^\ast\times\C\}$
which consists of $\{(a,0)\mid a\in\C^\ast\}$. We shall next find
an appropriate reduction of $\widetilde{\M}\rightarrow\M$ to an orbifold
principal $H$--bundle. We begin by giving a more detailed description of
the orbifold structure on $\M$ and the orbifold principal $G$--bundle
$\widetilde{\M}\rightarrow\M$.

First of all, we adopt the convention that $G$, as the automorphism
group of $\Sigma$, acts on $\Sigma$ from the left. Second, for the
orbifold structure on $\M$, we let $G$ act on $\widetilde{\M}$
from the left by defining $s\cdot f\equiv f\circ s^{-1}, \forall s\in G,
f\in\widetilde{\M}$. (This is because the convention is that the
group actions on a local uniformizing system are always from the
left.) To describe the orbifold structure, recall that for any
$f\in\widetilde{\M}$, there is a slice $S_f$ through $f$ which
has the following properties (cf \cite{Borel}):
\begin{itemize}
\item $S_f\subset\widetilde{\M}$ is a $2$--dimensional disc
containing $f$, which is invariant under the isotropy subgroup $G_f$
at $f$.
\item For any $s\in G$, $s\cdot S_f\cap S_f\neq\emptyset$ iff $s\in
G_f$.
\item There exists an open neighborhood $\O$ of $1\in G$ such that
the map $\phi_f\co \O\times S_f\rightarrow\widetilde{\M}$, defined by
$(s,h)\mapsto s\cdot h$, is an open embedding.
\end{itemize}
Let $U\equiv\bigsqcup_{f\in\widetilde{\M}} S_f$ be the disjoint union
of all slices. For any $h,h^\prime\in U$ which have the same
orbit in $\M$, and for any $s\in G$ such that $s\cdot h=h^\prime$,
let $\psi_{h^\prime,h}^s$ be the local self-diffeomorphism on $U$
defined as follows. Suppose $h\in S_f, h^\prime\in S_{f^\prime}$.
Then there is an open neighborhood $O_h\subset S_f$ of $h$, invariant
under the isotropy subgroup $G_h$ at $h$, such that $s\cdot O_h\subset
\phi_{f^\prime}(\O\times S_{f^\prime})$. Note that for any $g\in
O_h$, there is a unique $s^\prime\in\O$ and a unique $g^\prime\in
S_{f^\prime}$ such that $s\cdot g=\phi_{f^\prime}(s^\prime,g^\prime)
=s^\prime\cdot g^\prime$. We define $\psi_{h^\prime,h}^s(g)=g^\prime$,
which is clearly a local self-diffeomorphism on $U$ sending $h$ to
$h^\prime$. The orbifold structure on $\M$ is given by the pseudogroup
acting on $U$, which is generated by $\{\psi_{h^\prime,h}^s\}$.

To obtain the orbifold principal $G$--bundle $\widetilde{\M}\rightarrow\M$,
we let $G$ act on $\widetilde{\M}$ from the right by defining $f\cdot s
\equiv f\circ s,\forall s\in G, f\in\widetilde{\M}$. A local
trivialization of $\widetilde{\M}\rightarrow\M$ over a slice $S_f$ is
given by $(S_f\times G,G_f,\pi_f)$, where $G_f$ acts on $S_f\times G$
by $t\cdot (h,s)=(t\cdot h,ts)$, $\forall t\in G_f$, and where 
$\pi_f\co S_f\times G\rightarrow\widetilde{\M}$ sends $(h,s)$ to 
$h\cdot s=h\circ s$, which is invariant under the $G_f$--action
(note that $t\cdot h=h\circ t^{-1}$). The
transition function associated to each $\psi_{h^\prime,h}^s$ is given
by $g\mapsto\bar{\psi}_{h^\prime,h}^s(g)$, $\forall g\in
\mbox{Domain }(\psi_{h^\prime,h}^s)$, with $\bar{\psi}_{h^\prime,h}^s(g)\co 
G\rightarrow G$ being the multiplication by $(s^\prime)^{-1}s$ from left,
where $s^\prime\in\O$ is uniquely determined by $g\circ s^{-1}=
\psi_{h^\prime,h}^s(g)\circ (s^\prime)^{-1}$.

In the same vein, by letting $H$ act on $\widetilde{\M}$ from the
right, $\widetilde{\M}$ becomes an orbifold principal $H$--bundle
over $\widetilde{\M}\times_G (G/H)$. A reduction of
$\widetilde{\M}\rightarrow\M$ to an orbifold principal $H$--bundle is obtained
by taking a smooth section of $\widetilde{\M}\times_G (G/H)\rightarrow\M$.
Note that $G/H$ is naturally identified with $\C$, under which the
coset $(a,b)H$ goes to $b\in\C$. Now at any possible multiply
covered $f\in\widetilde{\M}$, $G_f\subset H$ iff $G_f$ is the cyclic
subgroup generated by $\mu_p$. Its action on $G/H$ is given by
$\mu_p\cdot (a,b)H=(\mu_p,0)(a,b)H$, which is simply the
multiplication by $\mu_p$ after identifying $G/H$ to $\C$.
Hence for any such $f$, a local uniformizing system of
$\widetilde{\M}\times_G (G/H)$ at $(f,0)$ is given by
$(S_f\times\C,G_f)$, where $G_f$ acts by $\mu_p\cdot (h,b)
=(\mu_p\cdot h,\mu_p b)$. To obtain a smooth section $u\co \M\rightarrow
\widetilde{\M}\times_G (G/H)$, we first pick a $G_f$--equivariant smooth
section $u_f\co S_f\rightarrow S_f\times\C$ for some arbitrary choice of a
multiply covered $f$ with $G_f\subset H$ (note that if there is such 
an $f$, its orbit in $\M$ is unique, cf Lemma 4.1 (1)), then extend 
it to the rest of $\M$, where $\widetilde{\M}\times_G (G/H)\rightarrow\M$ 
is an ordinary fiber bundle with a contractible fiber $\C$. We denote by
$\widehat{\M}\rightarrow\M$ the corresponding reduction to orbifold
principal $H$--bundle. Note that $\widehat{\M}$ is naturally
a 4--dimensional submanifold of $\widetilde{\M}$.

Fixing a choice of the reduction $\widehat{\M}\rightarrow\M$, we let
$Z\equiv\widehat{\M}\times_H\C$ be the associated orbifold complex
line bundle. Here $\C$ is canonically identified with $\Sigma\setminus
\{z_\infty\}$, and hence the action of $H$ on $\C$ is given by
complex multiplication.

There is a canonically defined smooth map of orbifolds
$\psi\co \widehat{\M}\times\Sigma\rightarrow X$, which
induces the evaluation map $(f,z)\mapsto f(z)$ between the
underlying spaces, cf Proposition 3.3.5 in Part I of \cite{C1}.
Note that each trivialization $S_f\times\C$ of $Z\rightarrow\M$
over a slice $S_f$ is a submanifold of $\widehat{\M}\times\Sigma$,
so that by restricting $\psi$ to $Z$, we obtain a smooth map of
orbifolds $\mbox{Ev}\co Z\rightarrow X$, which induces the evaluation
map $[(f,z)]\mapsto f(z)$ between the underlying spaces.

\begin{lem}
The map $\mbox{Ev}\co Z\rightarrow X$ is a diffeomorphism of orbifolds
onto $X\setminus\{x\}$.
\end{lem}

\proof
First of all, the map $\mbox{Ev}$ induces an injective map on the underlying
space. This is because each $J$--holomorphic curve parametrized by an
$f\in\widehat{\M}$ is a suborbifold, and any two distinct such $J$--holomorphic
curves $C,C^\prime$ intersect only at the singular point $x$. The latter
follows from the facts that (1) 
$C\cdot C^\prime\leq C_0\cdot C_0=\frac{p}{p+q}<1$,
so that by the intersection formula in Theorem 3.2, $C,C^\prime$ do
not intersect at any smooth point of $X$, (2) there is at most one such
$J$--holomorphic curve containing the other singular point $x^\prime$ of $X$.

Next we prove that the differential of $\mbox{Ev}$ is invertible
at each point of $Z$. Clearly the differential of $\mbox{Ev}$ is injective
along each fiber of $Z\rightarrow \M$, because each $f\in\widehat{\M}$ is
locally embedded on $\Sigma\setminus\{z_\infty\}$. Hence it suffices to
show that for any $f\in\widehat{\M}$ and any $u$ in the tangent space of
$\widehat{\M}$ at $f$ which is not tangent to the $H$--orbit through $f$,
$u(z)\in (TX)_{f(z)}$ is not tangent to $\mbox{Im }f$ for any
$z\in \Sigma\setminus\{z_\infty\}$. Note that $u$, being in the tangent
space of $\widehat{\M}$ at $f$, satisfies $D\underline{L}_f(u)=0$.

Now suppose to the contrary that $u$ is tangent to $\mbox{Im }f$
at some $z\in \Sigma\setminus\{z_\infty\}$. We can choose complex 
coordinates $w_1,w_2$ on a local uniformizing system at $f(z)$ 
such that $\mbox{Im }f$ is locally given by $w_2=0$, and $J$ equals 
the standard complex structure $J_0$ on $w_2=0$ (cf Lemma 1.2.2 in 
\cite{McD2}, or the corrected version of Lemma 2.5 in \cite{McD1}). 
Let $w=s+\sqrt{-1} t$ be a local holomorphic coordinate on $\Sigma$ 
centered at $z$, and set
$\partial=\frac{\partial}{\partial w}$, $\bar{\partial}=\frac{\partial}
{\partial \bar{w}}$.  Then
$$
\underline{L}(f)\equiv df+J\circ df\circ j=0, \;\forall f\in [\Sigma;X]
$$
can be written locally as
$$
\bar{\partial} f^i+a_{\bar{k}}^i(f)\bar{\partial}\bar{f}^k=0,
$$
where $f=(f^1,f^2)$, and $a_{\bar{k}}^i$ is a $2\times 2$ matrix of
smooth complex valued functions of $w_1,w_2$ which vanishes on $w_2=0$,
cf \cite{McD1}. Let $u_1,u_2$ be the components of $u$ in the
$\frac{\partial}{\partial w_1},\frac{\partial}{\partial w_2}$
directions, then $D\underline{L}_f(u)=0$ implies that
$$
\bar{\partial} u_2+ A u_2+B\bar{u}_2=0
$$
for some smooth complex valued functions $A,B$ of $s,t$. It follows easily
that $u_2$ satisfies
$$
|\Delta u_2|\leq c(|u_2|+|\partial_s u_2|+|\partial_t u_2|)
$$
pointwise for some constant $c>0$, where $\Delta=\partial^2_s+\partial^2_t$.
Note that $u_2$ is not constantly zero but $u_2(z)=0$ by the assumption,
hence by Hartman--Wintner's theorem \cite{HW},
$$
u_2(w)=aw^m+O(|w|^{m+1})
$$
for some nonzero $a\in\C$ and integer $m>0$.

Let $f_\lambda$, $\lambda\geq 0$, be a local smooth path in
$\widehat{\M}$ starting at $f$ which is tangent to $u$
at $\lambda=0$. Then in the local coordinate
system $\{w_1,w_2\}$, $f_\lambda$ is given by a pair of functions
$w_1=f^1_\lambda(w), w_2=f^2_\lambda(w)$ which satisfy
$$
(f^1_\lambda(w),f^2_\lambda(w))=\lambda (u_1(w),u_2(w))+ O(\lambda^2).
$$
We introduce $F_\lambda(w)=\lambda^{-1}(f^2_\lambda(w)-\lambda a w^m)$.
Then for any fixed, sufficiently small $\lambda\neq 0$, there is
an $r=r(\lambda)>0$ such that $|F_\lambda(w)|\leq |a|r^m$ for all $w$
satisfying $|w|\leq r$. For any such fixed $\lambda\neq 0$, we define 
a sequence $\{w=w_n\mid |w_n|\leq r=r(\lambda),\; n=1,2,\cdots\}$ 
inductively by solving
$$
F_\lambda(w_n)+aw_{n+1}^m=0,
$$
then $\{w_n\}$ has a limit $w_0$ in the disc $|w|\leq r=r(\lambda)$ 
satisfying
$$
F_\lambda(w_0)+aw_{0}^m=0.
$$
But this exactly means that $f^2_\lambda(w_0)=0$, which in turn
implies that $\mbox{Im }f_\lambda$ intersects with $\mbox{Im }f$
near $f(z)$, for any sufficiently small $\lambda\neq 0$. A contradiction.

Hence $u$ is nowhere tangent to $\mbox{Im }f$, and the differential
of $\mbox{Ev}\co Z\rightarrow X$ is injective, hence invertible by dimension
counting, at each point in $Z$.

To see that $\mbox{Ev}$ maps the underlying space of $Z$ onto that of
$X\setminus\{x\}$, note first that the image of $\mbox{Ev}$ is contained
in $X\setminus\{x\}$ and is an open subset. The latter is because the
differential of $\mbox{Ev}$ is invertible at each point of $Z$ so
that $\mbox{Ev}$ induces an open map between the underlying
spaces. On the other hand, the image of $\mbox{Ev}$ is also closed in
$X\setminus\{x\}$. To see this, suppose $\mbox{Ev}([(f_n,z_n)])=f_n(z_n)$
is a sequence of points in $X\setminus\{x\}$ which converges to $p\in
X\setminus\{x\}$. Since $\M$ is compact, a subsequence of $f_n$ (still
denoted by $f_n$ for simplicity) converges in $C^\infty$ to a
$f_0\in\widehat{\M}$ after reparametrization. If we let $z_0$ be
a limiting point of $z_n$ in $\Sigma$, then $z_0\neq z_\infty$, because
otherwise $p=\lim_{n\rightarrow\infty} f_n(z_n)=f_0(z_\infty)=x$, a
contradiction. This implies that the image of $\mbox{Ev}$ contains
$p=f_0(z_0)$, therefore it is closed in $X\setminus\{x\}$. Hence
$\mbox{Ev}$ maps $Z$ onto $X\setminus\{x\}$, and thus it
is a diffeomorphism from $Z$ onto $X\setminus\{x\}$.
\endproof

\sh{Proof of Theorem 1.1}

First of all, note that by Lemma 4.3, $\M$ is connected, and 
has an orbifold point of order $p$. The latter assertion is because 
there exists an $f\in\widehat{\M}$ such that $\mbox{Im }f$ contains 
the singular point $x^\prime\in X$, so that $f$ must be a multiply 
covered map. Moreover, $\M$ is orientable, and we shall orient $\M$
such that with the canonical orientation of orbifold complex line bundle
on $Z$, the map $\mbox{Ev}\co Z\rightarrow X$ is orientation-preserving. 
In order to determine the diffeomorphism type of $\M$ and the isomorphism 
class of the orbifold complex line bundle $Z\rightarrow\M$, we consider 
the family of regular neighborhoods of $x$:
$$
N_\epsilon\equiv\{(z_1,z_2)\mid |z_1|^2+|z_2|^2\leq \epsilon^2\}/\Z_{(p+q)}
$$
where $z_1,z_2$ are holomorphic coordinates on a local uniformizing
system at $x$ in which $C_0$ is locally given by $z_2=0$ and $C_0^\prime$,
the unique $J$--holomorphic curve containing both $x,x^\prime$, is locally
given by $z_1=0$.

\medskip

{\bf Claim}\qua {\sl
There exists an $\epsilon_0>0$ such that for any $0<\epsilon\leq\epsilon_0$,
$\partial N_\epsilon$ intersects transversely with each $J$--holomorphic 
curve in the family parametrized by $\M$ at a simple closed loop.
}

\proof 
For each $\lambda\in\M$, pick a local representative 
$(\hat{f}_\lambda,\rho_\lambda)$ of a member $f_\lambda\in\widetilde{\M}$
whose orbit in $\M$ is $\lambda$, and set $C_\lambda\equiv\mbox{Im }
f_\lambda$. Here $\rho_\lambda(\mu_{(p+q)})$
acts by $(z_1,z_2)\mapsto (\mu_{(p+q)}z_1,\mu_{(p+q)}^p z_2)$, and
$\hat{f}_\lambda=(U_\lambda,V_\lambda)$ for some holomorphic
functions $U_\lambda,V_\lambda$ defined on $D=\{z\in\C\mid |z|\leq 1\}$.
Observe (1) since $\M$ is compact, we may assume that for any sequence 
$\lambda_i\in\M$ converging to $\lambda_0\in\M$, there is a subsequence of 
$\lambda_i$, still denoted by $\lambda_i$, such that $\hat{f}_{\lambda_i}$ 
converges to $\hat{f}_{\lambda_0}\circ\xi$ for some holomorphic 
reparametrization $\xi$ of $D$, (2) for any $C_\lambda\neq C_0,C_0^\prime$, 
$C_\lambda\cdot C_0=\frac{p}{p+q}$ and $C_\lambda\cdot C_0^\prime
=\frac{1}{p+q}$, so that by the intersection formula in Theorem 3.2,
for any such a $\lambda$, $U_\lambda(z)=a_{\lambda,1}z+\cdots$, $V_\lambda(z)=
b_{\lambda,p}z^p+\cdots$ near $z=0$ for some $a_{\lambda,1}\neq 0$, 
$b_{\lambda,p}\neq 0$. (For $C_0$ or $C_0^\prime$, $\hat{f}_\lambda(z)$
equals $(a_{\lambda,1}z+\cdots,0)$ or $(0,b_{\lambda,p}z^p+\cdots)$ 
near $z=0$.)

Now for each $\lambda\in\M$, we write 
$U_\lambda(z)=a_{\lambda,1}z \cdot u_\lambda(z)$, 
$V_\lambda(z)=b_{\lambda,p}z^p\cdot v_\lambda(z)$ on $D$. Then there exist
$0<r_0\leq 1$, $0<\delta_0<1$, and $c>0$, which are independent of $\lambda$,
such that 
$$
1-\delta_0\leq|u_\lambda(z)|,|v_\lambda(z)|\leq 1+\delta_0, \mbox{ and }
|du_\lambda(z)|+|dv_\lambda(z)|\leq c
$$ 
when $|z|\leq r_0$. Write $z=r\exp(\sqrt{-1}\theta)$, and set 
$$
\mu_\lambda(r,\theta)\equiv |U_\lambda(z)|^2+|V_\lambda(z)|^2.
$$
Then each $\mu_\lambda$ is subharmonic on $D$, and a simple calculation 
shows that
$$
\frac{\partial\mu_\lambda(r,\theta)}{\partial r}
=|a_{\lambda,1}|^2 r(2|u_\lambda|^2+r\frac{\partial}{\partial r}|u_\lambda|^2)
+|b_{\lambda,p}|^2 r^{2p-1}(2p |v_\lambda|^2+r\frac{\partial}{\partial r}
|v_\lambda|^2),
$$
from which it follows that there exists $0<r_0^\prime\leq r_0$ such that
$$
\frac{\partial\mu_\lambda(r,\theta)}{\partial r}>0
$$ 
for all $\lambda\in\M$ whenever $0<r\leq r_0^\prime$.  

It remains to check that (1) there exists an $\epsilon_0>0$ such that
$\mu_\lambda(r,\theta)\leq \epsilon_0^2$ implies $r\leq r_0^\prime$, (2)
assuming the validity of (1), for any $0<\epsilon\leq\epsilon_0$,
the intersection of $\partial N_\epsilon$ with each $C_\lambda$, which is
transverse because $\frac{\partial\mu_\lambda}{\partial r}>0$ on
$\mu_\lambda^{-1}(\epsilon)$ by the validity of (1), is a simple closed 
loop. 

To see the former, note that $\mu_\lambda(r,\theta)\leq \epsilon_0^2$ 
implies 
$$      
r\leq \frac{2\epsilon_0}{|a_{\lambda,1}|\cdot |u_\lambda|+
(|b_{\lambda,p}|\cdot |v_\lambda|)^{1/p}},
$$
where on the other hand, it is easily seen that there exists a $c_1>0$
such that for any $\lambda\in\M$ and $|z|\leq r_0$,
$$
|a_{\lambda,1}|\cdot |u_\lambda|+(|b_{\lambda,p}|\cdot |v_\lambda|)^{1/p}
\geq c_1.
$$ 
To see the latter, suppose the intersection of $\partial N_\epsilon$ with 
some $C_\lambda$ consists of at least two components. Then either one of 
them bounds a disc in $D\setminus\{0\}$, or there is an annulus in 
$D\setminus\{0\}$ bounded by them. In any event, $\mu_\lambda$ will attain 
its minimum on the region at an interior point of the region (note that 
$\mu_\lambda$ is subharmonic on $D$), contradicting the fact that 
$\frac{\partial\mu_\lambda(r,\theta)}{\partial r}>0$ there. 
Hence the claim. 
\endproof

Back to the proof of Theorem 1.1. Let $E\rightarrow\M$ be 
the orbifold bundle of unit disc associated to $Z$. Then the
claim above implies that $X\setminus int(N_\epsilon)$ is diffeomorphic to
$E$ for any $0<\epsilon\leq\epsilon_0$. In particular, $\partial E$ is 
diffeomorphic to $\partial N_\epsilon=L(p+q,p)$. Note that $\partial E
\rightarrow \M$ defines a Seifert fibration of the lens space
$L(p+q,p)$ with one singular fiber of order $p$. Moreover, the Euler
number of the Seifert fibration, which equals the self-intersection of
the image of the zero section of $Z$ under the map $\mbox{Ev}\co Z\rightarrow
X$, is $1+\frac{q}{p}$ because it has a positive and transverse intersection 
with $C_0$ at a smooth point of $X$. This completely determines the 
diffeomorphism type of $\M$ and the isomorphism class of $Z$. 

Now observe that the same thing works for $X_{(p,q)}(R)$ as well.
In particular, the isomorphism class of $Z$ is independent of
$X$ and $X_{(p,q)}(R)$. Fix an $\epsilon>0$ and set $N\equiv N_\epsilon$.
Then from the proceeding paragraph, there are decompositions
$X=N\cup_{\phi_1} E$ and $X_{(p,q)}(R)=N\cup_{\phi_2} E$, where 
if we let $\gamma=\{z_2=0\}\cap\partial N$ and let $\gamma^\prime
=C_0\cap\partial E$, then $\phi_i(\gamma)=\gamma^\prime$, $i=1,2$.
Without loss of generality, we may assume $\phi_2=Id$ and $\gamma^\prime
=\gamma$ by fixing an identification of $\partial E$ with $\partial N$.
With this understood, we claim that $\phi_1$ is isotopic to the identity
through a family of diffeomorphisms $\phi_t\co \partial N\rightarrow
\partial N$ such that $\phi_t(\gamma)=\gamma$. 

First, assuming the validity of the claim, we obtain consequently
a diffeomorphism of orbifold pairs $\psi\co (X,C_0)\rightarrow 
(X_{(p,q)}(R),C_0)$, which preserves the singular point of order $p$ in 
$X$ and $X_{(p,q)}(R)$. By restricting $\psi$ to the complement of a regular 
neighborhood of the union of the singular point of order $p$ and the 
suborbifold $C_0$, we obtain a diffeomorphism $\psi^\prime\co W\rightarrow 
L(p,q)\times [0,1]$. 

It remains to verify the claim that $\phi_1$ is isotopic to the identity
through a family of diffeomorphisms $\phi_t\co \partial N\rightarrow
\partial N$ such that $\phi_t(\gamma)=\gamma$. To this end, let $Y$
be the complement of a regular neighborhood of $\gamma$ in $\partial N$.
Then $\pi_1(Y)$ is generated by the image of $\pi_1(\partial Y)$ in
$\pi_1(Y)$ induced by the inclusion $\partial Y\subset Y$, ie, $\pi_1(Y)$ 
is generated by the longitude and the meridian in $\partial Y\equiv T^2$. 
The diffeomorphism $\phi_1|_Y$ induces an automorphism of $\pi_1(Y)$ which 
is unique up to conjugation. In the present case, it is clear that the 
automorphism of $\pi_1(Y)$ can be chosen to be the identity map. Hence by the 
theorem of Waldhausen in \cite{Wal}, there exists an isotopy $\phi_t^\prime\co 
Y\rightarrow Y$ between $\phi_1|_Y$ and $Id$. Moreover, we may assume that 
$\phi_t^\prime|_{\partial Y}\co T^2\rightarrow T^2$ is given by a family of 
linear translations,
cf \cite{EE}. The latter implies particularly that $\phi_t^\prime$ can 
be extended to an isotopy $\phi_t$ from $\phi_1$ to $Id$ which satisfies 
$\phi_t(\gamma)=\gamma$. Hence the claim. 
\endproof

\sh{Proof of Corollary 1.2}

By Smith's theory (cf page 43 in \cite{Borel}), and by the assumption
that $\rho$ is free outside of a ball, we see easily that $\rho$ is
free in the complement of its fixed-point set, which consists of a
single point.  Then by applying (the proof of) Theorem 1.1 to the
quotient space of $\rho$, it follows easily that $\rho$ is conjugate
to a linear action by a diffeomorphism of $\R^4$. To see that the
diffeomorphism can be made identity outside of a ball, we note that in
the diffeomorphism $\psi\co (X,C_0)\rightarrow (X_{(p,q)}(R),C_0)$
constructed in the proof of Theorem 1.1, $\psi|_{C_0}\co
C_0\rightarrow C_0$ is isotopic to identity, from which it follows
easily.  \endproof

\end{document}